\documentclass[a4paper,12pt,twoside]{article}
\usepackage[left=2cm,right=2cm,top=3cm,bottom=3cm]{geometry}
\usepackage[utf8]{inputenc}
\usepackage[T1]{fontenc} 
\usepackage[english]{babel}
\usepackage{amsmath}
\usepackage{amsthm}
\usepackage{enumitem} 
\usepackage{stmaryrd} 
\usepackage{amssymb,bm} 
\usepackage{tikz,tkz-tab} 
\usepackage{multicol} 
\usepackage{parskip}
\usepackage[colorlinks=true, linkcolor=blue, citecolor=red]{hyperref} 
\usepackage{cleveref}
\usepackage{verbatim}
\usepackage{perpage} 
\usepackage[linesnumbered,ruled,vlined,french]{algorithm2e}
\usepackage{fancyhdr} 
\usepackage{titlesec}
\usepackage{setspace}
\usepackage{mathtools}
\usepackage{suffix}
\usepackage{tkz-euclide}
\usepackage{tikz}
	\usetikzlibrary{decorations.pathreplacing,decorations.markings, intersections, calc, angles, arrows}

\newcommand{\mycomment}[1]{}
\newcommand*\diff{\mathop{}\!\mathrm{d}}
\renewcommand{\mod}[1]{(\mathrm{mod}~#1)}
\WithSuffix\newcommand\psum*[1]{\mathop{\sum\nolimits^{\mathrlap{#1}}}}

\newcommand{\labo}{\vspace{1cm}
\textsc{Jérémy Dousselin : Université de Lorraine, CNRS, IECL, F-54000 Nancy, France}
\\ E-mail address : \texttt{jeremy.dousselin@univ-lorraine.fr}}

\newcommand{\CC}{\mathbb{C}}

\newcommand{\RR}{\mathbb{R}}

\newcommand{\I}{\mathcal{I}}

\newcommand{\meas}{\textrm{meas}}

\renewcommand{\S}{\mathcal{S}}

\newcommand{\C}{\mathcal{C}}

\newcommand{\E}{\mathcal{E}}

\newcommand{\M}{\mathcal{M}}

\renewcommand{\O}{\mathcal{O}}

\newcommand{\R}{\mathcal{R}}

\newcommand\numberthis{\addtocounter{equation}{1}\tag{\thesection.\theequation}}

\newcommand{\pb}{\overline{\O}}

\newcommand{\textAuthor}{To fill - Author}
\newcommand{\textTitle}{To fill - Title}
\renewcommand{\Res}{\text{Res}}
\newcommand{\AuthorIs}[1]{\renewcommand{\textAuthor}{#1}}
\newcommand{\TitleIs}[1]{\renewcommand{\textTitle}{#1}}
\newenvironment{entete}[1]
{
  \begin{center}
  {\Large\bf \uppercase{\textTitle}}
  \end{center}
  \vspace{0.5cm}
  \begin{center}
  {\large \uppercase{\textAuthor}}
  \end{center}
  \vspace{0.5cm}
	   {\begin{center}
	     \begin{minipage}{0.8\textwidth}
	     \small{ABSTRACT}.\footnotesize{~#1}\hfill
	     \end{minipage}
	     \end{center}
	      \vspace{0.2cm}
	      }
}

\renewcommand{\theequation}{\arabic{equation}}

\MakePerPage{footnote} 
\counterwithin*{equation}{section}
\labelformat{section}{Section #1}
\labelformat{subsection}{subsection #1}
\titleformat*{\section}{\large\bfseries\centering}
\newcounter{prop}[section]
\newcounter{theo}[section]

\newcounter{lemm}[section]
\newcounter{coro}[section]
\newcounter{fait}[section]
\newcounter{hypo}[section]
\newcounter{conj}[section]
\newcounter{glob}[section]
\newcounter{fig0}

\newenvironment{proo}[1]{\par\vspace{1em}\noindent  {\it Proof\ifx\relax#1\relax\else~(#1)\fi:}}{\qed\par}

\makeatletter
\newenvironment{fig}[1]{
	\refstepcounter{fig0}
	\protected@edef\@currentlabelname{Figure \arabic{fig0}}
	\newcommand{\copie}{#1}
	\begin{center}
}
{	
	\\
	Figure \arabic{fig0}: \copie
	\end{center}
}
\makeatother

\makeatletter

\newenvironment{theo}[1]{\par\vspace{1em}\noindent
\refstepcounter{glob}
\protected@edef\@currentlabelname{Theorem \arabic{section}.\arabic{glob}}\phantomsection
\noindent \textbf{Theorem \arabic{section}.\arabic{glob}}\ifx\relax#1\relax\else~(#1)\fi:\it}{
\par\vspace{1em}}

\newenvironment{coro}[1]{\par\vspace{1em} \noindent
\refstepcounter{glob}
\protected@edef\@currentlabelname{Corollary \arabic{section}.\arabic{glob}}\phantomsection
\noindent \textbf{Corollary \arabic{section}.\arabic{glob}}\ifx\relax#1\relax\else~(#1)\fi:\it}{\par\vspace{1em}
}

\newenvironment{lemm}[1]{ \par\vspace{1em} \noindent
\refstepcounter{glob}
\protected@edef\@currentlabelname{Lemma \arabic{section}.\arabic{glob}}\phantomsection
\noindent \textbf{Lemma \arabic{section}.\arabic{glob}}\ifx\relax#1\relax\else~(#1)\fi:\it}{\par\vspace{1em}
}

\newenvironment{rema}{\par\vspace{1em} \noindent
\refstepcounter{glob}
\protected@edef\@currentlabelname{Fact \arabic{section}.\arabic{glob}}\phantomsection\noindent\underline{\bf Remark \arabic{section}.\arabic{glob}}:}
	{\par\vspace{1em}
	}
\makeatother

\TitleIs{Zeros of linear combinations of Dirichlet $L$-functions on the critical line} 
\AuthorIs{Jérémy Dousselin} 

\begin{document}

\newcommand{\Ll}[2]{\log|L(1/2+i #1,\chi_{#2})|}

\entete{Let $F$ be a linear combination of $N\geq 1$ Dirichlet $L$-functions attached to even (or odd) primitive characters with the same modulus. Selberg proved that a positive proportion of non-trivial zeros of $F$ lie on the critical line. Our work here is to provide an explicit lower bound for this proportion. In particular, we show that the lower bound $2.16\times 10^{-6}/(N\log N)$ is admissible for large $N$.}

\section{Introduction}
\subsection{Zeros of the Riemann zeta-function on the critical line and the mollification method}
In 1838, Dirichlet introduced $L$-functions as a tool to study primes in arithmetic progression. Soon after, mathematicians understood how central these functions are in number theory, and especially how important it is to understand the distribution of their zeros. For Dirichlet $L$-functions, it is known that their non-trivial zeros lie in a strip, but a crucial hypothesis  states that all of these zeros should actually lie on the $\Re s=1/2$ line. Many ideas exist to study their zeros on this critical line, and one of them is to study the proportion of an $L$-function's zeros in the critical strip that are actually on the critical line. Showing that this proportion is non-zero is a thorny question answered by Selberg. The proof relies on a then recent method: the {\it mollification} method.

Let $N(T)$ be the number of zeros of the Riemann zeta-function $\zeta(s)$ that lie in the rectangle $\{s\in\CC:0\leq \Re(s)\leq 1,0\leq\Im(s)\leq T\}$, and let $N_0(T)$ be the number of these zeros on the critical line. In 1914, Hardy was the first to prove that infinitely many zeros of $\zeta(s)$ are on the critical line. Seven years later, he and Littlewood improved on this \cite{hlittlewood} and showed that we have 
\[N_0(T)\gg T\]
for all large $T$. However, since $N(T)\asymp T\log T$ (see Theorem 5.24 of \cite{iwaniec}), this lower bound is not strong enough to show that the proportion of non-trivial zeros lying on the critical line (we will call them {\it critical zeros}) is positive. It was Selberg \cite{selberg_zeta}, in 1942, who refined their argument by studying the zeros of $\zeta(s)|\eta(s)|^2$ instead of $\zeta(s)$, where $\eta(s)$ is a "mollifier". The role of this mollifying function is to prevent large values of $\zeta(s)$ from contributing too much, and hence making it possible to get sharper estimates. We will detail this mollifying process through the next section, when several notations are introduced. With this idea, originally due to Bohr and Landau \cite{blandau}, Selberg was able to prove that 
\[N_0(T)\gg T\log T,\]
which proves that $\kappa:=\liminf_T N_0(t)/N(T)>0$. His method actually produces a very small amount of zeros on the critical line, and hence yields a very small lower bound of $\kappa$. While we could not find the original proof, it is said (see page 68 of \cite{zhang} and $\mathsection 10.9$ of \cite{titchmarsh}) that Szu-Hoa Min computed the constant given by Selberg's method and proved that $\kappa\geq 1/60000$. This constant have been improved several times, still using the mollification method but not in the way Selberg did. The first breakthrough in this direction is due to Levison \cite{levinson} who proved that $\kappa\geq 1/3$, a result then refined by Conrey who showed that $\kappa\geq 2/5$.  Currently the record is hold by Pratt, Robles, Zaharescu and Zeindler who showed in 2019 that $\kappa\geq 5/12$ \cite{pratt}. In a few words, Levinson's and Selberg's methods are radically different. While Selberg's method detects very little zeros on the line, it is a safe method, for it detects zeros as the sign change of a real function: there is no risk of getting a negative lower bound. However, as explained through Section 3 of \cite{coniwa}, Levinson's method is more of a gamble. If the mollification is perfect, then one could potentially reach the $100\%$ lower bound. However, too crude estimates would lead to a negative bound for the counting of critical zeros. The counterpart of that risk is the reward, the method producing very high lower bound for the proportion of critical zeros. Note that the original proof of Levison is rather complicated, but a simpler one was given by Young in 2010 \cite{Young}.

In the case of a more general Dirichlet $L$-function both Selberg's method and Levinson's method still work, and the proof is not so different. Conrey, Iwaniec and Soundararajan \cite{consound} went a step further, and used Levison's method to prove that at least $56\%$ of zeros of the family of Dirichlet $L$-functions are on the critical line. This statement is to be understood as a sort of double average over both the $t$-aspect and the $q$-aspect of a Dirichlet $L$-function, mostly focusing on the latter. 

\subsection{The case of linear combinations of Dirichlet $L$-functions}
The study of zeros of linear combinations of $L$-functions is motivated by the existence of certain zeta functions (called Epstein zeta functions) which satisfy all properties of $L$-functions, except the existence of an Euler product representation. It turns out that these zeta functions possess non-trivial zeros off the critical line and hence do not satisfy the Generalized Riemann Hypothesis. However, these zeta functions can be expressed as linear combinations of Hecke $L$-functions, which is why several number theorists investigated the zeros of general linear combinations of $L$-functions.

Montgomery conjectured that under natural conditions (notably a certain notion of independence of these $L$-functions), 100\% of the zeros of a linear combination of $L$-functions lie on the critical line (which implies that the counterexamples to the Generalized Riemann Hypothesis are rare). In their celebrated paper, Bombieri and Hejhal \cite{bombieri} proved -under some reasonable, but yet to be proven, hypotheses- that this is true.

Unconditionally, Karatsuba \cite{kara} tried to tackle the question and proved that for a specific linear combination of two $L$-functions, we have $N_0(T)\gg T(\log T)^{1/2-\varepsilon}$, which narrowly fails to show that the proportion of critical zeros is positive. Surprisingly, a few years later Selberg proved \cite{selberg},\cite{sel_comp} that his method actually applies in the case of "any" linear combination of Dirichlet $L$-functions. Recently in 2016, Rezvyakova \cite{irina_deg2} proved that Selberg's method also applies in the case of a linear combination of $L$-functions of degree two attached to automorphic forms.

Our goal here is to provide an explicit lower bound for the proportion of zeros on the critical line for $F$, a linear combination of Dirichlet $L$-functions attached to even (or odd) primitive characters with the same modulus. Apart from trivial zeros implied by the functional equation, the zeros of $F$ lie in a vertical strip. Indeed, $F(s)$ has a series representation for $\Re(s)>1$. Since the first non-zero term dominates the others, $F(s)$ has no zero if $\Re(s)$ is large enough, and hence we may define 
\[\sigma_F:=\sup\{\Re(s):F(s)=0\}.\]
From the functional equation satisfied by $F(s)$, we may deduce the set of zeros on $\Re(s)<1-\sigma_F$, called {\it trivial zeros}. The other zeros are confined in the strip $1-\sigma_F\leq \Re(s)\leq \sigma_F$, and we call them {\it non-trivial zeros}. We denote by $N(T,F)$ the number non-trivial zeros with imaginary part in $(0,T)$, and by $N_0(T,F)$ the number of these zeros that are on the critical line. We also write 
\[\kappa_F:=\liminf_{T}\frac{N_0(2T,F)-N_0(T,F)}{N(2T,F)-N(T,F)}.\]
In order to give an unconditional lower bound for $\kappa_F$, we will use the mollification method {\it à la} Selberg to prove the following:
\begin{theo}{}\label{mainth}
Let $F$ be a linear combination of $N\geq 1$ distinct Dirichlet $L$-functions attached to even primitive characters with the same modulus. Then for any $\frak A>1/\varkappa$, any $\varkappa\in(0,1/8)$ and any $\theta\in (0,1)$:
\[\kappa_F\geq 2\pi\left(\frac{1}{2\frak A}-4N\frac{\frak C_1(\frak A)+\frak C_2 }{\frak A^3}\right).\]
Here, $\frak C_1(\frak A)$ and $\frak C_2$ are defined in \nameref{prop:esti_inte} and implicitly depend on $\theta$ and $\varkappa$. If $N=1$, this can be improved by \eqref{borne_meilleure}.
\end{theo}
From this, some computations and optimisations lead to the following, more explicit theorem.
\begin{theo}{}\label{application}
Let $F$ be a linear combination of $N\geq 1$ distinct Dirichlet $L$-functions attached to even primitive characters with the same modulus. For small $N$'s, we get the following lower bounds for $\kappa_F$, as well as the corresponding $\frak A$ and $\theta$.
\vspace{0.3cm}\\
\begin{center}
\begin{tabular}{|c|c|c|c|}
\hline $N$ & $\frak A$ & $\theta$ & $\kappa_F\geq $
\\\hline $1$ & $29056699.107509706$ & $0.011$ & $5.45\times 10^{-8}$
\\\hline $2$ & $212583177.09901848$ & $0.0016$ & $7.38\times 10^{-9}$
\\\hline $3$ & $319102776.4709714$ & $0.0014$ & $4.91\times 10^{-9}$
\\\hline $4$ & $425715589.6389222$ & $0.0013$ & $3.68\times 10^{-9}$
\\\hline $5$ & $532459869.61320543$ & $0.0012$ & $2.94\times 10^{-9}$
\\\hline $10$ & $1067086846.4520979$ & $0.001$ & $1.46\times 10^{-9}$
\\\hline $100$ & $10776391786.558016$ & $0.0004$ & $1.45\times 10^{-10}$
\\\hline $1000$ & $109024453631.91109$ & $0.0002$ & $1.43\times 10^{-11}$
\\\hline
\end{tabular}
\end{center}
\vspace{0.3cm}
Moreover, when we are interested in very large $N\geq 3$, we may state the following: for any given $0<\varepsilon<1/3$, we have
\[N\geq\frac{2.9\times 10^{-11}}{\varepsilon^3} \implies \kappa_F\geq \frac{2.161\times 10^{-6}}{N\log N}\left(1-3\varepsilon-\frac{1.14\log\log N}{\log N}-\frac{862}{\log N}-\frac{1.2\times 10^7\varepsilon}{\log^2N}\right).\]
\end{theo}
Note that the case of odd characters is handled similarly, and we only deal with the even ones for the sake of simplicity.
\begin{rema}
Selberg \cite{selberg},\cite{sel_comp} only proved that $\kappa_F\geq c/N^2$, for some $c>0$. However, he mentioned that he could strengthen the lower bound to $c/(N\log N)$, but he did not provide any idea on how one can achieve such a result. By getting a better error term that him in a crucial lemma (compare our \nameref{intracine} to the Lemma 9 of \cite{selberg_zeta}), we were able to retrieve the result. We will give more technical details later in the proof, but essentially Selberg's estimate corresponding to our lemma was given with an error term of $\O(\sqrt \frak A)$, while our error term is $\O(1)$. This made him find $\frak C_1(\frak A)\asymp \frak A^{3/2}$ instead of our $\frak C_1(\frak A)\asymp\frak A\log\frak A$, which explains the improved lower bound.
\end{rema}
\begin{rema}
In fact, we will prove a slightly more general result, but chose to present a simpler statement for simplicity. The two theorems above actually hold for any linear combination of the form \eqref{ex_F} below, thus removing the need for the characters to share the same modulus. In either the simpler or the more general statement, the important property that $F$ must satisfy is a functional equation, which is inherited from the $L$-functions' functional equation.
\end{rema}
\noindent\underline{\bf Acknowledgement:} I would like to thank my PhD advisor Youness Lamzouri for suggesting this problem to me, for his time and guidance, and for several improvement suggestions. Un remerciement à Pascal Ciot pour une aide bibliographique bien utile.
\section{Notation and outline of the proof}\label{sec-notation}
Now and for the rest of the paper, we write $f(T)\lesssim g(T)$ to say that $f(T)\leq g(T)(1+o(1))$. Similarly, $f(T)\gtrsim g(T)$ means that $g(T)\lesssim f(T)$. The $o$ symbol is to be understood as valid when $T\to\infty$. Several times throughout the paper, the constant implied by the $\O$-symbols may depend on the conductors of the considered $L$-functions, and we will never explicit this dependence. We will also write $f=\pb(g)$ if $f\lesssim g$. We define $\log^+(x)=\max(0,\log x)$.

Let $N\geq 1$ be an integer, $\chi_1,...,\chi_N$ be even distinct primitive Dirichlet characters to moduli $q_1,...,q_N$ respectively. We know (see section 5.1 of \cite{iwaniec}) that there is a complex number $\varepsilon_j$ with $|\varepsilon_j|=1$ such that 
\[\phi(s,\chi_j):=\varepsilon_j\pi^{-s/2}q_j^{s/2}\Gamma\left(\frac{s}{2}\right)L(s,\chi_j)\]
satisfies $\phi(s,\chi_j)=\overline{\phi(1-\overline s,\chi_j)}$ for all $1\leq j\leq N$. Write 
\[\numberthis\label{ex_F} F(s):=\sum_{j=1}^N\frac{c_j}{q_j^{1/4}}\varepsilon_jq_j^{s/2} L(s,\chi_j),\]
where $c_j\in\RR^*$\footnote{We chose to normalize each summand by $q_j^{1/4}$ because $q_j^{s/2-1/4}$ is real on the critical line. This has no impact on the proof since $c_j$ are arbitrarily chosen.}. Standard arguments show (see Theorem 5.8 of \cite{iwaniec} for the case of a single $L$-function) that 
\[\numberthis\label{eq:N}N(2T,F)-N(T,F)=\frac{T}{2\pi}\left(\log T+B_F\right)+\O(\log T),\]
for some $B_F\in\RR$.

The idea is to count zeros of $F$ on the critical line thanks to sign changes of a good function. Let $1\leq j\leq N$, and write 
\[\vartheta(s):=\arg \pi^{-s/2}\Gamma\left(\frac s2\right)\]
and
\[X_j(s):=\varepsilon_j q_j^{s/2-1/4}e^{i\vartheta(s)}L(s,\chi_j).\]
This function $X_j$, thanks to the functional equation given by $\phi$, takes real values on the critical line. Moreover, $1/2+it$ is a zero of $X_j$ if and only if $1/2+it$ is a zero of $L$. Thus, any sign change of $X_j(1/2+i\cdot)$ implies the existence of a zero of $L$ on the critical line.

Let $T$ be a large real number, and we define $(\alpha_j(n))$ by
\[L(s,\chi_j)^{-1/2}=\sum_n\frac{\alpha_j(n)}{n^s},\text{ for }\Re(s)>1.\]
For $T\leq t\leq 2T$, $\xi=T^{\varkappa}$ with $\varkappa>0$ to be chosen later\footnote{One may keep in mind that $\varkappa$ can be taken as close as we want to $1/8$, as we will show later.}, we write
\[\eta_j(s):=\sum_{n\leq \xi}\frac{\alpha_j(n)}{n^{s}}\M(n)=:\sum_{n}\frac{\beta_j(n)}{n^{s}},\]
where $\mathcal M$ is the continuous function
\[\M(x):=\mathcal M(\xi,x):=\frac{1}{\log \xi}\left\{\begin{array}{cl}
\log \xi&\text{ if }1\leq x\leq \xi^\theta,
\\\frac{1}{1-\theta}\log(\xi/x)&\text{ if } \xi^\theta< x\leq \xi,
\\ 0&\text{ if }x>\xi,
\end{array}\right.\]
where $0<\theta<1$ is a parameter to be chosen later. Also note that $|\beta_j(n)|\leq |\alpha_j(n)|\leq 1$ since $\alpha_j(n)=\tau_{-1/2}(n)\chi_j(n)$. Here $\tau_z$ is the $z$-th divisor function, defined to be the coefficient of the Dirichlet series of $\zeta^z$, satisfying the inequality $|\tau_z(n)|\leq \tau_{|z|}(n)$.

Selberg worked with a slightly different function in \cite{selberg_zeta}, choosing the continuous function
\[\M_{Sel}(\xi,x):=\frac{1}{\log \xi}\left\{\begin{array}{cl}
\log(\xi/x)&\text{ if } 1\leq x\leq \xi,
\\ 0&\text{ if }x>\xi.
\end{array}\right.
\]
This will not heavily impact the rest for the proof, for we have 
\[\numberthis\label{lien-molli}\M(\xi,x)=\frac{1}{1-\theta}\left(\M_{Sel}(\xi,x)-\theta\M_{Sel}(\xi^\theta, x)\right).\]

Selberg \cite{sel_comp} suggested to use this function in the case $\theta=1/2$ instead of $\M_{sel}$. This is also done in \cite{hb},\cite{irina_deg2}. The pros of working with this function are that the upper bound \eqref{estiM} below does not depends on $N$, making it better than the original one found by Selberg in \cite{selberg_zeta}, when $N$ is large. Furthermore, the proof of \eqref{estiM} is also greatly simplified. The drawback of using $\M$ instead of $\M_{sel}$ is that the constant term $\frak C_2$ is rather large, which worsen our lower bound for small $N$'s. See for example the case $N=1$, where the lower bound $\kappa_\zeta \geq 1/60000$ was found admissible, while our lower bound for $N=1$ is roughly $10^{-7}$. However, these results for $N=1$ are both vastly improved by Levinson's, and hence the real interest of our Theorem is when $N>1$.

One of our improvements is to introduce the parameter $\theta$, which has no impact on the complexity of the proof but has a strong impact on the quality of our lower bound. As the computations of \nameref{application} show, it seems that the optimal $\theta$ is relatively small, and that it gets smaller as $N$ gets larger.

Now, we introduce three integrals that we will use to detect sign changes of $X_j$ on the critical line. Let $H=\frak A/\log T$, where $\frak A> 1/\varkappa$ is a quantity depending only on $N$, $\varkappa$ and $\theta$ to be chosen later. We define 
\[I_j(t,H):=\int_t^{t+H}X_j(1/2+iu)|\eta_j(1/2+u)|^2\diff u,\]
\[J_j(t,H):=\int_t^{t+H}|X_j(1/2+iu)\eta^2_j(1/2+u)|\diff u,\]
and 
\[M_j(t,H):=\int_t^{t+H}L(1/2+iu,\chi_j)\eta^2_j(1/2+iu)\diff u-H.\]
It is clear that 
\[\numberthis\label{min-J}J_j(t,H)\geq H-|M_j(t,H)|,\]
thus if 
\[\numberthis\label{cond-suff}|M_j(t,H)|+|I_j(t,H)|<H,\]
we get that 
\[\numberthis\label{cond-end}J_j(t,H)>|I_j(t,H)|,\]
which immediately implies that there is a sign change of $X_j(1/2+iu)$ in $(t,t+H)$, and hence a zero of $L(1/2+it,\chi_j)$ in $(t,t+H)$.

This detection method works well on short intervals at one condition: the function must not have very large spikes. Indeed, if the curve of $X_j$ were composed of one large spike and many smaller spikes, then the area under this curve would be concentrated under the larger spike, making all the remaining area almost undetectable by integral computations, therefore making $I_j$ and $J_j$ very close even though there could be many sign changes. That is why we mollify the function by a positive factor, so that we get a new function with the same sign changes but much less large oscillations. Without mollification, we would still be able to detect some sign changes, but not accurately enough on such small (of length $H\asymp 1/\log T$) intervals. \nameref{fig:moli} illustrates this phenomenon by comparing the function $\zeta(s)$, the function $\zeta(s)$ mollified by $\eta(s)$, and the function $\zeta(s)$ mollified by Selberg's mollifier.
\begin{fig}{Mollification at work on $\zeta(s)$. Here we chose the simplest $\theta=1/2$ and $\varkappa=1/8$.}\label{fig:moli}
\includegraphics[scale=0.5]{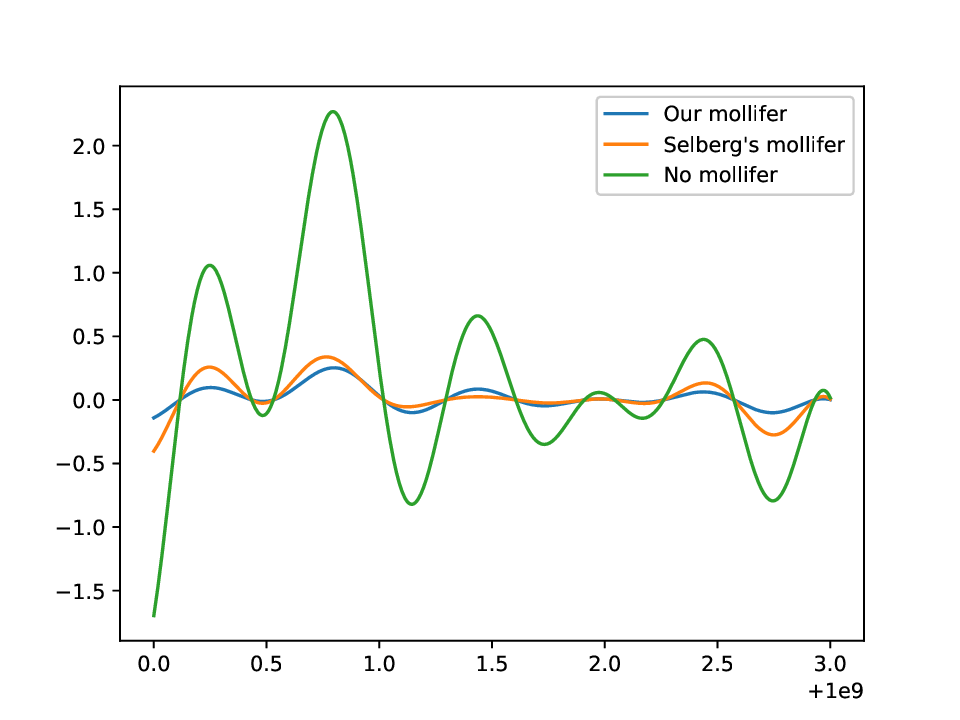}
\end{fig}\vspace{0.5cm}

Now to exploit \eqref{cond-suff} and \eqref{cond-end}, we rely on the following estimates, which are the main results of this paper.
\begin{theo}{}\label{prop:esti_inte}
Suppose that $\varkappa<1/8$ and $0<\theta<1$. Recall that $H=\frak A/\log T$. Then we have 
\[\numberthis\label{estiI}\int_T^{2T}|I_j(t,H)|^2\diff t\lesssim \frak C_1(\frak A)\frac{T}{(\log T)^2},\]
and 
\[\numberthis\label{estiM}\int_T^{2T}|M_j(t,H)|^2\diff t\lesssim \frak C_2\frac T{(\log T)^2}.\]
Here,
\[\numberthis\label{eq:defC3}\frak C_3=\left(\frac{1}{8\varkappa}+\frac{3}2\right)\left(\frac{e^{\varrho_\theta}+e^{\varrho_\theta\theta}}{(1-\theta)\sqrt {\varrho_\theta\pi}}\frac{\Gamma(1/4)}{\Gamma(3/4)}\right)^4\prod_{p}\left(1+\frac{3p^2-3p+1}{p^4-3p^3+3p^2-p}\right),\]
with $\varrho_\theta$ being the only positive solution of $-1+2\theta x+e^{x(1-\theta)}(-1+2x)=0$, and
\[\frak C_2=\frac{6(\frak C_3+1+2\sqrt{\frak C_3} )}{(\theta\varkappa)^2}.\]
Furthermore,
\[\frak C_1(\frak A)=8\frak C_5^2\left(\frak K_1 \frak A\log\frak A+\frak K_2\frak A+\frak K_3\log \frak A+\frak K_4\right),\]
where $\frak K$'s are given in \eqref{deffrakK} and under. 
\end{theo}
\begin{rema}
Rather than focusing on the complex and lengthy explicit expression for $\frak C_1$, one should note that $\frak C_1$ is around $\frak A\log\frak A$ (multiplicative constant aside) when $\frak A$ is large. This is the important piece of information conveyed by the expression above. Also note that the estimate \eqref{estiM} holds whenever $\varkappa<2/9$.
\end{rema}

To illustrate how the method works, we will apply it to a single $L$-function called $F$. Since we are working with a single $L$-function, we drop the indices $j$ for simplicity, and we will do it several times later without mentioning this again. Let $\S$ be the subset of $(T,2T)$ consisting of $t$'s such that $|I(t,H)|=J(t,H)$. Then
\[\int_\S |I(t,H)|\diff t=\int_\S J(t,H)\diff t,\]
and hence, by \eqref{estiI} and the Cauchy-Schwarz inequality
\[\int_\S J(t,H)\diff t\lesssim \left(\meas(\S)\frak C_1(\frak A)\frac{T}{(\log T)^2}\right)^{1/2},\]
where $\meas$ is Lebesgue's measure on $\RR$. Moreover, by \eqref{min-J} and \eqref{estiM}, we find that 
\[\int_\S J(t,H)\diff t\gtrsim \meas(\S) H-\left(\meas(\S)\frak C_2\frac T{(\log T)^2}\right)^{1/2}.\]
Thus by combining these we find 
\[\meas(\S)\lesssim \frac{(\sqrt{\frak C_1(\frak A)}+\sqrt{\frak C_2})^2}{\frak A^2}T.\]
Now we divide $(T,2T)$ into $\lfloor T/(2H)\rfloor$ pairs of abutting intervals $\I_1,\I_2$ of length $H$ each (except maybe for the last $\I_2$). For each pair, there is at least one critical zero of $F$ in $\I_1$ or $\I_2$ unless $\I_1\subset \S$. If $n_1$ is the number of these bad intervals $\I_1$, then we find 
\[n_1\times H\leq \meas(\S)\lesssim\frac{(\sqrt{\frak C_1(\frak A)}+\sqrt{\frak C_2})^2}{\frak A^2}T ,\]
and hence 
\[n_1\lesssim \frac{(\sqrt{\frak C_1(\frak A)}+\sqrt{\frak C_2})^2}{\frak A^3}T\log T.\]
Therefore there are at least 
\[\left\lfloor\frac{T}{2H}\right\rfloor -n_1\gtrsim \left(\frac{1}{2\frak A}-\frac{(\sqrt{\frak C_1(\frak A)}+\sqrt{\frak C_2})^2}{\frak A^3}\right)T\log T\]
zeros of $F$ on the critical line in $(T,2T)$. By \eqref{eq:N}, this implies that
\[\numberthis\label{borne_meilleure}\kappa_F\geq 2\pi\left(\frac{1}{2\frak A}-\frac{(\sqrt{\frak C_1(\frak A)}+\sqrt{\frak C_2})^2}{\frak A^3}\right).\]
Optimising the parameters $\theta$, $\varkappa$ and $\frak A$ would imply the result in the case of a single $L$-function. For a general linear combination of $L$-functions, the idea is to prove that the interval $(T,2T)$ can be split into "good" subsets where one of the $L$-functions dominates the others, and "bad" subsets whose measure we can control. Then, on the good subsets we can work as if we only had a single $L$-function. We will explain this with further details in the next section.

\begin{rema}
This lower bound is slightly stronger than the general one presented in \nameref{mainth}, and the loss comes from a slight technical difficulty arising in the case of a combination of multiple $L$-functions. 
\end{rema} 
\section{How to handle the general case}
We first suppose that \nameref{prop:esti_inte} is true and we detail how to use it to prove \nameref{mainth}. We follow the main lines of approach of Selberg \cite{sel_comp}. First, we need the following version of Selberg's central limit theorem for a difference of two distinct $L$-functions.
\begin{theo}{}\label{th:cent}
Let $\chi_1$ and $\chi_2$ be distinct primitive Dirichlet characters. Let $a<b$ be two real numbers and $\mathfrak 1_{[a,b]}$ be the characteristic function of $[a,b]$. Then as $T\to\infty$ we have
\[\int_T^{2T}\mathfrak 1_{[a,b]}\left(\frac{\Ll t1-\Ll t2}{\sqrt{2\pi \log\log T}}\right)\diff t= T\int_a^b e^{-\pi u^2}\diff u+\O\left(T\frac{\log\log\log^2T}{\sqrt{\log\log T}}\right).\]
\end{theo}
\begin{proo}{}
Tsang proved a similar result for the Riemann zeta-function in his thesis \cite{tsang}. Very few modifications are needed to adapt the proof to the case of a difference of two Dirichlet $L$-functions.  
\end{proo}
We will also need the following lemmas.
\begin{lemm}{}\label{lemma-tsang}
Fix $\chi$ an even Dirichlet character, let $k\geq 1$ be fixed, and let $x=T^{1/(100k)}$. We have that
\[\int_T^{2T}\left|\Ll t{}-\Re\sum_{p<x}\frac{\chi(p)}{p^{1/2+it}}\right|^{2k}\diff t\ll TB^kk^{4k},\]
for some constant $B>0$.
\end{lemm}
\begin{proo}{}
Again, Tsang proved it for the Riemann zeta-function in his thesis (see \cite{tsang}, Theorem 5.1). The proof remains unchanged for a Dirichlet $L$-function.
\end{proo}
\begin{lemm}{Lemma 3 of \cite{sound}}\label{lem-sound}
Let $T$ be large and let $2\leq x\leq T$. Let $k$ be a natural number such that $x^k\leq T/\log T$. For any complex numbers $a(p)$ we have
\[\int_T^{2T}\left|\sum_{p\leq x}\frac{a(p)}{p^{1/2+it}}\right|^{2k}\diff t\ll Tk!\left(\sum_{p\leq x}\frac{|a(p)|^2}p\right)^k.\]
\end{lemm}
\begin{lemm}{}\label{lemma-stretch}
Fix $1\leq j\leq N$ and $k\geq 1$. Define, for $t\in (T,2T)$
\[\Delta_j(t,H):=\frac{1}{H}\int_t^{t+H}\Ll uj\diff u.\]
Then, uniformly for $0\leq v\leq H$, we have that 
\[\int_{T}^{2T}(\Delta_j(t,H)-\Ll{(t+v)}j)^{2k}\diff t\ll TB^kk^{4k},\]
for some constant $B>0$ that depends on $\frak A$.
\end{lemm}
\begin{proo}{}
Put $x=T^{1/(100k)}$. Define
\[Err(t):=\Ll tj-\Re\sum_{p<x}\frac{\chi(p)}{p^{1/2+it}}.\]
We have
\begin{align*}
&\int_{T}^{2T}(\Delta_j(t,H)-\Ll{(t+v)}j)^{2k}\diff t\ll 
\\&\int_T^{2T}\left[\frac{1}{H}\int_t^{t+H}Err(u)\diff u\right]^{2k}\diff t+\int_T^{2T}\left[\frac{1}H\int_t^{t+H}\Re\sum_{p<x}\frac{\chi(p)}{p^{1/2+iu}}\diff u-\Re\sum_{p<x}\frac{\chi(p)}{p^{1/2+i(t+v)}}\right]^{2k}\diff t
\\&\numberthis\label{eq:intDelta}+\int_T^{2T}Err(t+v)^{2k}\diff t.
\end{align*}
Now remember that $H\leq 1$ since $T$ is large. By Hölder's inequality and switching the integrals, the first integral is 
\[\ll\frac{1}H\int_T^{2T}\int_t^{t+H}Err(u)^{2k}\diff u\diff t\ll\int_T^{2T+1}Err(u)^{2k}\diff u.\]
Thus, \nameref{lemma-tsang} shows that both the first and the third integrals are 
\[\numberthis\label{eqcentral}\ll TB^kk^{4k},\]
for some constant $B>0$.

Now we have to deal with the second integral of \eqref{eq:intDelta}, that we write as 
\[\int_2:=\int_T^{2T}\left[\Re\sum_{p<x}\frac{\chi(p)\left(\frac 1H\int_0^{H}p^{-iu}\diff u-p^{-iv}\right)}{p^{1/2+it}}\right]^{2k}\diff t.\]
By \nameref{lem-sound}, we know that we only have to keep the diagonal terms of this:
\[\numberthis\label{eq:int2}\int_2\ll Tk!\left(\sum_{p<x}\frac{1}p\left|\frac{p^{-iH}-1}{-iH\log p}-p^{-iv}\right|^2\right)^k.\]
We split the inner sum in two parts, the first over $p<e^{1/H}$ and the second over the remaining $p$'s. By Mertens' Theorem  and Taylor expanding $p^{-iH}$ and $p^{-iv}$, we find that the contribution of the first part is 
\begin{align*}
&\ll \sum_{p<e^{1/H}}\frac1p\left|\frac{-iH\log p+\O(H^2\log^2p)}{-iH\log p}-1+\O(H\log p)\right|^2
\\&\ll \sum_{p<e^{1/H}}\frac 1p\left|H\log p\right|^2\ll H\sum_{p<e^{1/H}}\frac{\log p}p\ll 1.
\end{align*}
Again by Mertens' Theorem, the sum over larger $p$'s is 
\[\ll \sum_{e^{1/H}\leq p<x}\frac{1}{p}\ll \log\left(H\log x\right)\ll \log(H\log T)=\log \frak A.\]
Putting this back in \eqref{eq:int2} and \eqref{eq:intDelta} together with \eqref{eqcentral}, we prove the lemma.
\end{proo}
\begin{lemm}{}\label{lemme-mollifer}
We suppose that $\varkappa<2/9$ and $0<\theta<1$. Then we have, for $1\leq j\leq N$,
\[\int_T^{2T}|L(1/2+it,\chi_j)\eta_j^2(1/2+it)|^2\diff t\lesssim \frak C_3 T\]
where $\frak C_3$ is defined in \eqref{eq:defC3}.
\end{lemm}
The quite technical proof of this lemma is postponed until the end of the paper, for it contains some arguments also found in the proof of \eqref{estiI}. We can now proceed to the proof of \nameref{mainth}.

\begin{proo}{of \nameref{mainth}}
Let $\delta>0$ be small ($\delta=1/10$ should be enough), and fix $1\leq j\neq j'\leq N$. Let $E_{j,j'}$ be the set consisting of $t\in(T,2T)$ such that 
\[|\Ll tj-\Ll t{j'}|\leq (\log\log T)^\delta.\]
If $\epsilon=(\log\log T)^\delta/\sqrt{2\pi\log\log T}$, then \nameref{th:cent} implies that  
\begin{align*}
\meas(E_{j,j'})&=\int_T^{2T}\mathfrak 1_{[-\epsilon,\epsilon]}\left(\frac{\Ll tj-\Ll t{j'}}{\sqrt{2\pi \log\log T}}\right)\diff t
\\&=T\int_{-\epsilon}^{\epsilon} e^{-\pi u^2}\diff u+\O\left(T\frac{\log\log\log^2T}{\sqrt{\log\log T}}\right)\ll T (\log \log T)^{\delta-1/2}.
\end{align*}

Outside the union of these exceptional sets, $(T,2T)$ may be split into subsets on which only one single $L$-function decisively dominates the others. Of course, for each said subset the dominating $L$-function may vary. Now we shall see that this dominance is stable over relatively long (compared to $1/\log T$) stretches.

We integrate \nameref{lemma-stretch} over $0\leq v\leq H$, and we get that 
\[\numberthis\label{moment}\int_{T}^{2T}\int_0^H(\Delta_j(t,H)-\Ll{(t+v)}j)^{2k}\diff v\diff t\ll HTB^kk^{4k},\]
for some constant $B>0$ and any fixed integer $k\geq 1$. Now write $W_j(t)$ the subset of $v\in[0,H]$ such that
\[\numberthis\label{eq_deltaj}|\Delta_j(t,H)-\Ll{(t+v)}j|>(\log\log T)^{\delta/2}.\]
Taking $k=\lfloor 6/\delta\rfloor+1$ so that $k\delta\geq 6$, \eqref{moment} implies that 
\[\numberthis\label{mesureW}\meas(W_j(t))\leq \frac{H}{(\log\log T)^3}\]
except maybe for a subset of $t\in(T,2T)$ of measure 
\begin{align*}
&\leq \int_{T}^{2T}\meas(W_j(t))\frac{(\log \log T)^3}H\diff t
\\&\leq \int_T^{2T}\int_0^H |\Delta_j(t,H)-\Ll{(t+v)}j|^{2k}\frac{(\log \log T)^{3-k\delta}}H\diff v\diff t\ll \frac{T}{(\log\log T)^3}.
\end{align*}

Finally by \nameref{lemma-stretch} with $v=0$, it follows that the set $F_j$ consisting of $t\in(T,2T)$ such that 
\[\numberthis\label{eq_delta_F_stretch} |\Delta_j(t,H)-\Ll tj|>(\log\log T)^{\delta/2}\]
has measure 
\begin{align*}
\meas(F_j)\leq \int_{T}^{2T}\left(\frac{|\Delta_j(t,H)-\Ll tj|}{(\log\log T)^{\delta/2}}\right)^{2k}\diff t\ll \frac{T}{(\log\log T)^6}.
\end{align*}

Now if we exclude from $(T,2T)$ all the $t$ such that $t\in E_{j,j'}$ for some $1\leq j\neq j'\leq N$, or such that $t\in F_j$ for some $j$, or such that $\meas(W_j(t))>H/(\log\log T)^3$ for some $j$, we get that $(T,2T)$ except for a subset of size $\O(T/(\log\log T)^{1/2-\delta})$ can be divided into $N$ subsets $S_j$ such that for all $t\in S_j$ and all $j'\neq j$:
\[\numberthis\label{eq_LL1}\Ll tj-\Ll t{j'}> (\log\log T)^\delta,\]
and such that for any $t\in S_j$ and any $u\in H_t:=(0,H)-\bigcup_{r=1}^NW_r(t)$:
\begin{align*}
&\Ll {(t+u)}j-\Ll{(t+u)}{j'}= \Bigl[\Ll{(t+u)}j-\Delta_j(t, H)\Bigr]
\\&-\Bigl[\Ll{(t+u)}{j'}-\Delta_{j'}(t,H)\Bigr]+\Bigl[\Delta_j(t,H)-\Ll tj\Bigr]
\\&+\Bigl[\Ll t{j'}-\Delta_{j'}(t,H)\Bigr]+\Ll tj-\Ll t{j'}
\\\numberthis\label{eq:int_Ht}&>(\log\log T)^\delta-4(\log\log T)^{\delta/2}>\frac12(\log\log T)^\delta.
\end{align*}
The first two differences are lower bounded by \eqref{eq_deltaj}, the next two are bounded by \eqref{eq_delta_F_stretch}, while the last one is bounded by \eqref{eq_LL1}. 

From \nameref{lemme-mollifer}, we see that 
\[\numberthis\label{intmollif}\int_t^{t+H}|L(1/2+iu,\chi_j)\eta_j^2(1/2+iu)|^2\diff u\leq H\log\log T\]
except for a subset of $t\in(T,2T)$ of measure
\begin{align*}
&\leq \frac{1}{H\log\log T}\int_T^{2T}\int_t^{t+H}|L(1/2+iu,\chi_j)\eta_j^2(1/2+iu)|^2\diff u\diff t
\\&\ll\frac{1}{\log\log T}\int_T^{2T+1}|L(1/2+iu,\chi_j)\eta_j^2(1/2+iu)|^2\diff u\ll \frac{T}{\log \log T}.
\end{align*}
Let us denote by $S_j^*$ the set $S_j$ deprived of those $t$'s. Note that by construction, we have that 
\begin{align*}
 \sum_{m=1}^N\meas(S_m^*)&=\sum_{m=1}^N\meas(S_m)+\O\left(\frac{T}{\log\log T}\right)
\\&=T+\O\left(\frac{T}{(\log\log T)^{1/2-\delta}}+\frac{T}{\log\log T}\right)
\\\numberthis\label{sizeS}&\sim T.
\end{align*}

Let $\mathcal F:t\mapsto\sum_{n=1}^Nc_nX_n(1/2+it)$, and for $t\in S_j^*$ we define 
\[I_j^*(t,H):=\int_{H_t}\mathcal F(t+u)|\eta_j^2(1/2+i(t+u))|\diff u\]
and
\[J_j^*(t,H):=\int_{H_t}\left|\mathcal F(t+u)\eta_j^2(1/2+i(t+u))\right|\diff u.\]
Any sign change of $\mathcal F$ would imply a zero of $e^{i\vartheta(1/2+i\cdot)}F(1/2+i\cdot)$, and hence a critical zero of $F$. Again, we will use these integrals above to detect sign changes of $\mathcal F$. We have, for $t\in S_j^*$,
\begin{align*}
&I_j^*(t,H)=c_j\int_{H_t}X_j(1/2+i(t+u))|\eta_j^2(1/2+i(t+u))|\diff u
\\&\qquad+(N-1)\O\left(\int_{H_t}|L(1/2+i(t+u),\chi_j)\eta_j^2(1/2+i(t+u))|e^{-(\log\log T)^\delta/2}\diff u\right)
\\&=c_j\int_{0}^HX_j(1/2+i(t+u))|\eta_j^2(1/2+i(t+u))|\diff u
\\&\qquad+\int_{\cup_{r=1}^NW_r(t)}X_j(1/2+i(t+u))|\eta_j^2(1/2+i(t+u))|\diff u
\\&\numberthis\label{eq:ierr}\qquad+\O\left(e^{-(\log\log T)^\delta/2}\sqrt{\int_t^{t+H}|L(1/2+iu,\chi_j)\eta_j^2(1/2+iu)|^2\diff u}\sqrt{\meas(H_t)}\right).
\end{align*}
The first equality relies on \eqref{eq:int_Ht}, while the second equality comes from an application of the Cauchy-Schwarz inequality. By the Cauchy-Schwarz inequality, \eqref{mesureW} and \eqref{intmollif}, the second integral of this is
\begin{align*}
\ll \sqrt{\int_t^{t+H}|L(1/2+iu,\chi_j)\eta_j^2(1/2+iu)|^2\diff u}\sqrt{\sum_{i=1}^N\meas(W_i(t))}\ll\frac{H}{\log\log T}.
\end{align*}
The same bound holds for the error term of \eqref{eq:ierr}. Therefore,
\[\numberthis\label{I*}I_j^*(t,H)=c_jI_j(t,H)+\O\left(\frac{H}{\log\log T}\right).\]
Similarly for $t\in S_j^*$,
\[\numberthis\label{J*}J_j^*(t,H)=|c_j|J_j(t,H)+\O\left(\frac{H}{\log\log T}\right).\]

Now we need a slight change compared to the case of a single $L$-function. Indeed, we are now working with sets $S_j^*$ that are probably not intervals, and this prevents us from using the exact method that we described in \ref{sec-notation}. However, we may easily overcome this difficulty. We fix $\varepsilon>0$. For any $t\in (T,2T)$ except for a subset $\E_j$ of measure $\lesssim (2+\varepsilon)^2(\frak C_1(\frak A)+\frak C_2)T/(H\log T)^2=(2+\varepsilon)^2(\frak C_1(\frak A)+\frak C_2) T/\frak A^2$, \nameref{prop:esti_inte} implies that 
\[|I_j(t,H)|\leq \frac{H}{2+\varepsilon},\text{ and }|M_j(t,H)|\leq \frac{H}{2+\varepsilon}.\]
Thus, for any $t\in S_j^*\setminus\E_j$ we have by \eqref{I*}:
\begin{align*}
\frac{|I_j^*(t,H)|}{|c_j|}&= |I_j(t,H)|+\O\left(\frac H{\log\log T}\right)\lesssim\frac {H}{2+\varepsilon}.
\end{align*}
On the other hand by \eqref{min-J} and \eqref{J*}, 
\[\frac{J_j^*(t,H)}{|c_j|}=J_j(t,H)+\O\left(\frac{H}{\log\log T}\right)\geq H-|M_j(t,H)|+\O\left(\frac{H}{\log\log T}\right)\gtrsim H-\frac{H}{2+\varepsilon}.\]
Combining these two inequalities above, we find that for any $t\in S_j^*\setminus\E_j$, $F(1/2+i\cdot)$ has a zero in $(t,t+H)$, since $|I^*_j(t,H)|<J_j^*(t,H)$ if $T$ is large enough. Summing over $j$ and using \eqref{sizeS}, we find that for any $t$ in a subset of $(T,2T)$ of measure 
\[\gtrsim\sum_{j}\meas(S_j^*)-N(2+\varepsilon)^2\frac{(\frak C_1(\frak A)+\frak C_2) T}{\frak A^2}\gtrsim \left(1-N(2+\varepsilon)^2\frac{\frak C_1(\frak A)+\frak C_2 }{\frak A^2}\right)T\]
$F(1/2+i\cdot)$ has a sign change in $(t,t+H)$. Again we divide $(T,2T)$ into $\lfloor T/(2H)\rfloor$ pairs of abutting intervals $\I_1,\I_2$ of length $H$ each (except maybe for the last $\I_2$). For each pair, there is at least one zero in $\I_1$ or $\I_2$ unless $\I_1$ is included in a subset $\E_\infty$ of measure $\lesssim N(2+\varepsilon)^2\frac{\frak C_1(\frak A)+\frak C_2 }{\frak A^2}T$. If $n_1$ is the number of these $\I_1$, then we find 
\[n_1\times H\lesssim \meas(\E_\infty)\lesssim N(2+\varepsilon)^2\frac{\frak C_1(\frak A)+\frak C_2 }{\frak A^2}T ,\]
and hence 
\[n_1\lesssim N(2+\varepsilon)^2\frac{\frak C_1(\frak A)+\frak C_2 }{\frak A^3}T\log T.\]
Therefore there are at least 
 \[\left\lfloor\frac{T}{2H}\right\rfloor -n_1\gtrsim \left(\frac{1}{2\frak A}-N(2+\varepsilon)^2\frac{\frak C_1(\frak A)+\frak C_2 }{\frak A^3}\right)T\log T\]
zeros of $F$ on the critical line in $(T,2T)$. $\varepsilon$ being as small as we want, we may use \eqref{eq:N} to conclude.
\end{proo}
\section{Preliminary results}
We first need a few lemmas to deal with the final proofs. For the remaining proofs, recall that $(\beta(n))$ is only supported on integers $n\leq \xi$. Thus, we will often not specify the bounds in our summations indexes to lighten notations. Several lemmas among the followings are refinements of lemmas found in \cite{selberg_zeta}.
\begin{lemm}{}\label{lemmeA}
For integers $a,b\geq 1$, we put $a|b^\infty$ if $p|a\implies p|b$ for all primes $p$. Let $d\leq T^{2\varkappa}$ and $q\geq 1$ be positive integers, $(q,d)=1$, and let $\delta|d^\infty$. For $0\leq v \leq H$ a real number, we define
\[S_d(\delta,v):=\sum_{(n,qd)=1}\frac{\mathcal M(n\delta)\tau_{-1/2}(n\delta)}{n^{1+iv}}.\]
Then for any $c>0$,
\begin{align*}
|S_d(\delta,v)|\lesssim \frac{\tau_{1/2}(\delta)}{1-\theta}\frac{e^{c}+e^{c \theta}}{2\sqrt{\pi\varkappa c}}\left(\sqrt{\frac{v\log T}{c}}\sqrt{\pi\varkappa}+\frac{\Gamma(1/4)}{\Gamma(3/4)}\right)\sqrt{\frac{dq}{\phi(d)\phi(q)\log T}}.
\end{align*}
In particular this holds for the optimal $c=\rho(v\log T,\theta)$, the only positive solution in $X$ of $e^{(1-\theta)X}\left[2X(a+b\sqrt X)-2a-b\sqrt X\right]+2\theta X(a+b\sqrt X)-2a-b\sqrt X=0$. Here, 
\[a:=\sqrt{\pi\varkappa v\log T },\qquad b:=\frac{\Gamma(1/4)}{\Gamma(3/4)}.\]
\end{lemm}
\begin{proo}{}
We know that for any $\varepsilon>0$,
\[\frac{1}{2\pi i}\int_{\varepsilon-i\infty}^{\varepsilon+i\infty}\frac{x^s}{s^2}\diff s=\left\{\begin{array}{cl}
0&\text{if }0<x\leq 1,
\\ \log x& \text{if }x>1.
\end{array}\right.\]
Therefore, by \eqref{lien-molli}, our function $\mathcal{M}$ can be expressed as the following Mellin representation
\[\mathcal{M}(x)=\frac{1}{2\pi i}\int_{\Re(s)=\varepsilon}x^{-s}\frac{1}{(1-\theta)\log \xi} \frac{\xi^s-\xi^{s\theta}}{s^2}\diff s.\]
Since for any integer $n\geq 1$ we have $(n,d)=1\implies (n,\delta)=1$, we deduce that
\begin{align*}
S_d(\delta,v)&=\frac{1}{2\pi i(1-\theta)(\log \xi)}\sum_{(n,qd)=1}\frac{\tau_{-1/2}(n\delta)}{n^{1+iv}}\int_{\Re(s)=\varepsilon}(n\delta)^{-s}\frac{\xi^s-\xi^{s\theta}}{s^2}\diff s
\\\numberthis\label{eqad1}&=\frac{\tau_{-1/2}(\delta)}{2\pi i(1-\theta)\varkappa(\log T)}\int_{\Re(s)=\varepsilon}Z(1+iv+s)\frac{\xi^s-\xi^{s\theta}}{\delta^s s^2}\diff s,
\end{align*}
where 
\[Z(s)=\sum_{(n,qd)=1}\frac{\tau_{-1/2}(n)}{n^s}=\prod_{p\nmid qd}(1-p^{-s})^{1/2}=\prod_{p| qd}(1-p^{-s})^{-1/2}\zeta(s)^{-1/2}.\]
We take $\varepsilon=c/\log T$, for some $c>0$ to be chosen later. We want to remove the factor $\zeta(s)$ by taking advantage of the estimate $\zeta(1+s)\approx 1/s$ when $s$ is small. Therefore, we will show that the \eqref{eqad1} integral can be truncated to only keep relatively small imaginary part.

Since $1/\zeta(1+s+iv)\ll |s|+|v|$ and on the line $\Re(s)=\varepsilon$,
\[\numberthis\label{inegProd}\left|\prod_{p| qd}(1-p^{-(s+1)})^{-1/2}\right|\leq \prod_{p| qd}(1-p^{-1})^{-1/2}=\sqrt{\frac{dq}{\phi(dq)}}=\sqrt{\frac{dq}{\phi(d)\phi(q)}},\]
we may bound the contribution of the integral of \eqref{eqad1} on $|\Im(s)|> 1/\log\log T$ by
\[\numberthis\label{int_petit_im}\O\left(\frac{\tau_{1/2}(\delta)}{\log T}\sqrt{\frac{d}{\phi(d)}}\int_{\substack{\Re(s)=\varepsilon\\|\Im(s)|> 1/\log\log T}}\frac{\sqrt{|s|+|v|}|\diff s|}{|s|^{2}}\right)=\O\left(\tau_{1/2}(\delta)\sqrt{\frac{d}{\phi(d)}}\frac{\sqrt{\log\log T}}{\log T}\right).\]
Moreover, we have 
\begin{align*}
&\int_{\substack{\Re(s)=\varepsilon\\|\Im(s)|\leq 1/\log\log T}}\frac{(|s|+v)^{3/2}}{|s|^2}|\diff s|\ll \int_{\substack{\Re(s)=\varepsilon\\|\Im(s)|\leq 1/\log\log T}}\frac{|s|^{3/2}+v^{3/2}}{|s|^2}|\diff s|
\\&\ll \sqrt\varepsilon\int_{t\leq 1/(\varepsilon\log\log T)}\frac{1}{(1+t^2)^{1/4}}\diff t+\frac{H^{3/2}}{\varepsilon}\int_{t\leq 1/(\varepsilon\log\log T)}\frac{1}{1+t^2}\diff t\ll\frac{1}{\sqrt{\log\log T}}+\frac{1}{\sqrt{\log T}}
\\&\ll\numberthis\label{intermediaire} \frac{1}{\sqrt{\log\log T}}
\end{align*}
Using the fact that $1/\zeta(1+z)^{1/2}=\sqrt{z+\O(z^2)}=\sqrt z+\O(|z|^{3/2})$ when $z\to 0$, $\Re(z)>0$, we may use \eqref{inegProd}, \eqref{int_petit_im} and \eqref{intermediaire} in \eqref{eqad1} to get 
\begin{align*}
S_d(\delta,v)&=\frac{\tau_{-1/2}(\delta)}{1-\theta}\frac{1}{\varkappa(\log T)2\pi i}\int_{\substack{\Re(s)=\varepsilon\\|\Im(s)|\leq 1/\log\log T}}\prod_{p| qd}(1-p^{-(1+iv+s)})^{-1/2}\sqrt{iv+s}\frac{\xi^s-\xi^{s\theta}}{\delta^s s^2}\diff s
\\&\numberthis\label{plusdenom}+\O\left(\tau_{1/2}(\delta)\sqrt{\frac{d}{\phi(d)}}\frac{\sqrt{\log\log T}}{\log T}\right).
\end{align*} 
We recall that $c=\varepsilon\log T$. Using \eqref{inegProd}, we see that the absolute value of the main term of this is
\begin{align*}
&\lesssim \frac{\tau_{1/2}(\delta)}{2\pi \varkappa (1-\theta)}\sqrt{\frac{qd}{\phi(d)\phi(q)}}(e^{c\varkappa}+e^{\varkappa c\theta})\frac{1}{\sqrt{\log T}}\int_{\substack{\Re(s)=\varepsilon\\|\Im(s)|\leq 1/\log\log T}}\frac{\sqrt v+\sqrt{|s|}}{|s|^2}|\diff s|
\\&\lesssim \frac{\tau_{1/2}(\delta)}{2\pi \varkappa (1-\theta)}\sqrt{\frac{qd}{\phi(d)\phi(q)}}(e^{c\varkappa}+e^{\varkappa c\theta})\frac{1}{\log T}\left(\frac{\sqrt{v}}\varepsilon\int_\RR\frac{\diff t}{1+t^2}+\frac{1}{\sqrt \varepsilon}\int_\RR \frac{\diff t}{(1+t^2)^{3/4}}\right).
\end{align*}
Of course, we have $\int_\RR\diff t/(1+t^2)=\pi$. We let $B(x,y)=\int_0^{\infty}\frac{s^{y-1}}{(1+s)^{x+y}}\diff s$ be the Beta-function, $x,y>0$. By property of this Beta-function (see Theorem II.0.8 of \cite{tenenbaum}), we deduce that
\[\frac{\sqrt{\pi}\Gamma(1/4)}{\Gamma(3/4)}=B(1/4,1/2)=\int_0^{\infty}\frac{s^{-1/2}}{(1+s)^{3/4}}\diff s=\int_0^{\infty}\frac{2}{(1+t^2)^{3/4}}\diff t=\int_\RR\frac{1}{(1+t^2)^{3/4}}\diff t.\]
Combining these and going back to \eqref{plusdenom},
\[|S_d(\delta,v)|\lesssim \frac{\tau_{1/2}(\delta)}{1-\theta}\frac{e^{\varkappa c}+e^{\varkappa c \theta}}{2\pi \varkappa\sqrt c}\left(\sqrt{\frac{v\log T}{c}}\pi+\sqrt\pi\frac{\Gamma(1/4)}{\Gamma(3/4)}\right)\sqrt{\frac{dq}{\phi(d)\phi(q)\log T}}.\]
By choosing the optimal $c$, one concludes.
\end{proo}
\begin{lemm}{}\label{intracine}
Fix $0< v\leq H$. Suppose that $\xi/\delta \geq 2$, where $\delta>0$ is an integer. Put $\varepsilon=c/\log T$, where $c> 0$. Then for $r=1,2,3$: 
\[\frac{1}{2\pi i}\int_{\substack{\Re(s)=\varepsilon\\|\Im(s)|\leq 1}}\sqrt{s+iv}\left(\frac{\xi}{\delta}\right)^s\frac{1}{s^r}\diff s=-\frac{v^{3/2-r}}{2\sqrt{\pi}}\Delta_r(\log(\xi/\delta)v)+\O\left(\frac{1}{\log(\xi/\delta)}\right),\]
where 
\[\Delta_1(X):=-2X^{-1/2}+\int_0^{X} \frac{e^{-it}-1}{t^{3/2}}\diff t,\]
and for $r=2,3$:
\[\Delta_r(X):=\int_0^{X}(X-x)^{r-2}\left(-2x^{-1/2}+\int_0^x \frac{e^{-it}-1}{t^{3/2}}\diff t\right)\diff x.\]
\end{lemm}
\begin{proo}{} Let $L$ be a curve starting at $-\infty-i$, joining the $\Re(s)=\varepsilon$ axis horizontally, going up to $\varepsilon+i $ as a semi-circle and then returning back to $-\infty +i$. \vspace{0.5cm}

\begin{fig}{The continuous line is $L$, and the dotted one is Hankel's contour.}
\begin{tikzpicture}[scale=0.5]
\tkzDefPoints{0/0/O,1.4/0.5/A, 1.4/-0.5/B, 3.7/-0.5/C, 3.7/0.5/D,2/0/E, 2/2/F, 2/-2/G, 2/0/H, 0/2/I, 0/-2/J}
\tkzDefPoints{2/0/OO,-1.38/0.5/AAA, -1.38/-0.5/BBB, -7/-0.5/CCC, -7/0.5/DDD, -7/2/FFF, -7/-2/GGG}
\tkzDefPoints{0/0/O,-1.2/0.5/AA, -1.2/-0.5/BB, -7/-0.5/CC, -7/0.5/DD}

\tkzDrawCircle[draw=gray,densely dotted, line width=1pt, postaction={decorate,decoration={
        markings,
        mark=at position 0.1 with {\arrow{stealth}}
      }}](O,A)
\tkzDrawPolygon[draw=white, fill=white](AA,BB,CC,DD)
\draw[draw=gray, densely dotted, line width=1pt, postaction={decorate,decoration={
        markings,
        mark=at position .5 with {\arrow{stealth}}
      }}] (AAA) -- (DDD);
\draw[draw=gray, densely dotted, line width=1pt, postaction={decorate,decoration={
        markings,
        mark=at position .5 with {\arrow{stealth}}
      }}] (CCC) -- (BBB);
\draw[->] (-5,0) -- (5,0);
\draw[->] (0,-5) -- (0,5);
\tkzDefPointsBy[rotation=center O angle 120](A){GA}
\draw[gray, ->] (O) -- (GA) node[midway, above]{\tiny $r$};
\tkzLabelPoint[above right](O){\small $0$}
\tkzLabelPoint[above right](I){\small $i$}
\tkzLabelPoint[above right](H){\small $\varepsilon$}
\tkzDrawPoint(H)
\tkzLabelPoint[below right](J){\small $-i$}
\draw[draw=black, line width=1pt, postaction={decorate,decoration={
        markings,
        mark=at position .5 with {\arrow{stealth}}
      }}] (F) -- (FFF);
\draw[draw=black, line width=1pt, postaction={decorate,decoration={
        markings,
        mark=at position .5 with {\arrow{stealth}}
      }}] (GGG) -- (G);
\begin{scope}
    \clip (2,-3) rectangle (5,3);
    \tkzDrawCircle[draw=black, line width=1pt, postaction={decorate,decoration={
        markings,
        mark=at position 0.1 with {\arrow{stealth}}
      }}](OO,G);
\end{scope}
\end{tikzpicture}
\end{fig}
\vspace{0.5cm}
Note that for any $s\in L$, $|s|>v$. For any integer $j$ and any complex $z$, we recall the definition $\binom z j:=\frac{z(z-1)(z-2)...(z-j+1)}{j!}=(-1)^{j-1}\frac{z(1-z)(2-z)...(j-1-z)}{j!}$. Then for any $s\in L$,
\[\sqrt{iv+s}=\sqrt s\left(1+\frac{iv}s\right)^{1/2}=\sum_{j=0}^\infty\frac{(iv)^j}{s^{j-1/2}}\binom{1/2}j.\]
We chose that branch of $\sqrt{s+iv}$ which is positive if $iv+s$ is real and positive. Observe that for any integer $\ell\geq 1$: 
\[\Gamma(\ell-1/2)=\left(\ell-\frac32\right)\left(\ell-\frac52\right)\times \cdots\times\frac32\times\frac12\Gamma(1/2)=2\Gamma(1/2) \ell! (-1)^{\ell-1}\binom{1/2}\ell.\]
Therefore, for any $j\geq 1$, we find that
\begin{align*}
\Gamma(j+r-1/2)&=\left(j+r-\frac32\right)\left( j+r-\frac52\right)\times\cdots\times \left(j-\frac12\right)\Gamma(j-1/2)
\\&=2\Gamma(1/2) j! (-1)^{j-1}\binom{1/2}j\left(j+r-\frac32\right)\left( j+r-\frac52\right)\times\cdots\times \left(j-\frac12\right).
\end{align*}
This holds if $j=0$. Let $\log 2\leq u\leq \log \xi$. Then by Hankel's formula for $1/\Gamma$ together with the above formulas, we find that
\begin{align*}
\frac{1}{2\pi i}\int_L\sqrt{iv+s}\frac{e^{us}}{s^r}\diff s&=\sum_{j=0}^\infty(iv)^j\binom{1/2}j\frac{1}{2\pi i}\int_L\frac{e^{us}}{s^{j+r-1/2}}\diff s
\\&=\sum_{j=0}^\infty(iv)^j\binom{1/2}j\frac{u^{j+r-3/2}}{\Gamma(j+r-1/2)}
\\&\numberthis\label{eqchemin}=-\frac{u^{r-3/2}}{2\Gamma(1/2)}\sum_{j=0}^\infty\frac{(-iuv)^j}{j!(j-1/2)(j+1/2)...(j+r-3/2)}.
\end{align*}
By Taylor expanding the exponential below, we find that 
\begin{align*}
(uv)^{1/2}\left(-2(uv)^{-1/2}+\int_0^{uv} \frac{e^{-it}-1}{t^{3/2}}\diff t\right)&=(uv)^{1/2}\left(-2(uv)^{-1/2}+\sum_{j=1}^\infty \frac{(-i)^j(uv)^{j-1/2}}{j!(j-1/2)}\right)
\\&=\sum_{j=0}^\infty\frac{(-iuv)^j}{j!(j-1/2)},
\end{align*}
while for $r>1$, Taylor expanding the exponential below together with $r-2$ successive integrations by parts lead to
\begin{align*}
&\frac{(uv)^{3/2-r}}{(r-2)!}\int_0^{uv}(uv-x)^{r-2}\left(-2x^{-1/2}+\int_0^x \frac{e^{-it}-1}{t^{3/2}}\diff t\right)\diff x
\\&\qquad\qquad=\frac{(uv)^{3/2-r}}{(r-2)!}\int_0^{uv}(uv-x)^{r-2}\sum_{j=0}^\infty\frac{(-i)^jx^{j-1/2}}{j!(j-1/2)}\diff x
\\&\qquad\qquad=(uv)^{3/2-r}\int_0^{uv}\sum_{j=0}^\infty\frac{(-i)^jx^{j+r-5/2}}{j!(j-1/2)(j+1/2)...(j+r-5/2)}.
\\&\qquad\qquad=\sum_{j=0}^\infty\frac{(-iuv)^j}{j!(j-1/2)(j+1/2)...(j+r-3/2)}
\end{align*}

Using \eqref{eqchemin} and deforming the semi-circle part of $L$ as the segment $[\varepsilon-i,\varepsilon+i]$, one finds 
\[\numberthis\label{sansnom}\frac{1}{2\pi i}\int_{\substack{\Re(s)=\varepsilon\\|\Im(s)|\leq 1}}\sqrt{iv+s}\frac{e^{us}}{s^r}\diff s=-\frac{v^{3/2-r}}{2\sqrt\pi}\Delta_r(uv)+\O\left(\int_{\substack{\Re(s)\leq\varepsilon\\ \Im(s)=1}}\frac{\sqrt v+\sqrt{|s|}}{|s|^r}e^{u\Re(s)}|\diff s|\right).\]
Since we assumed that $u\leq \log \xi \ll 1/\varepsilon$ and since $r\geq 1$, the error term is 
\[\ll \int_{-\infty}^\varepsilon \frac{1+(1+t^2)^{1/4}}{(1+t^2)^{r/2}}e^{ut}\diff t\ll \frac{e^{u\varepsilon}}u\ll \frac1u.\]
Using \eqref{sansnom} with $u=\log(\xi/\delta)$, one concludes.
\end{proo}
\begin{lemm}{}\label{lemm_esti_delta}
Let $X$ be a positive real number. We have:
\[\Delta_2(X)=\left\{\begin{array}{ll}
-4X^{1/2}+\pb\left(\frac43 X^{3/2}\right) & \text{if }X\leq 1,
\\ 2\sqrt\pi e^{-3i\pi /4}(X-1)-4+\pb\left(\frac{16}3\right) & \text{else.}
\end{array}\right.\]
Thus for any $X\geq 0$, 
\[\Delta_2(X)=2\sqrt\pi e^{-3i\pi/4}(X-1)-4+\pb\left(\frac{16}3\right).\]
Furthermore if $X\leq 1$,
\[\Delta_3(X)=-\frac{8}3X^{3/2}+\pb\left(\frac8{15} X^{5/2}\right).\]
\end{lemm}
\begin{proo}{}
If $x\leq 1$ then the inequality $|e^{-it}-1|\leq t$, valid for $t\geq 0$, implies 
\[\int_0^x \frac{e^{-it}-1}{t^{3/2}}\diff t=\pb(2\sqrt x).\]
If $x>1$ then an integration by parts leads to
\begin{align*}
-2x^{-1/2}+\int_0^x \frac{e^{-it}-1}{t^{3/2}}\diff t&=\int_0^\infty\frac{e^{-it}-1}{t^{3/2}}\diff t-2x^{-1/2}-\int_{x}^\infty\frac{e^{-it}-1}{t^{3/2}}\diff t
\\&=\int_0^\infty\frac{e^{-it}-1}{t^{3/2}}\diff t+\pb(2x^{-3/2}).
\end{align*}
Observe that an integration by parts and a change of variable lead to
\[\int_0^\infty\frac{e^{-it}-1}{t^{3/2}}\diff t=-2i\int_0^\infty\frac{ e^{-it}}{\sqrt t}\diff t=-4i\int_0^\infty e^{-iu^2}\diff u.\]
The last integral happens to be Fresnel's integral, whose value is known to be $e^{-i\pi/4}\frac{\sqrt \pi}2$. Therefore,
\[\int_0^\infty\frac{e^{-it}-1}{t^{3/2}}\diff t=2\sqrt{\pi}e^{-3i\pi / 4}.\]
Thus, if $0\leq X\leq 1$,
\[\Delta_2(X)=\int_0^X\left(-2x^{-1/2}+\pb(2\sqrt x)\right)\diff x=-4\sqrt X+\pb\left(\frac 43X^{3/2}\right),\]
and if $X\geq 1$:
\[\Delta_2(X)=-4+\pb\left(\frac43\right)+\int_1^X\left(2\sqrt \pi e^{-3i\pi/4}+\pb(2x^{-3/2})\right)\diff x=2\sqrt \pi e^{-3i\pi/4}(X-1)-4+\pb\left(\frac43\right)+\pb(4).\]
This proves the first part of the lemma.

For the second part of the lemma, we know that for any $X\geq 1$,
\[\Delta_2(X)=2\sqrt\pi e^{-3i\pi/4}(X-1)-4+\pb\left(\frac{16}3\right).\]
For $X\leq 1$, we have 
\begin{align*}
\Delta_2(X)&=-4\sqrt X+\pb(4/3)
\\&=2\sqrt\pi e^{-3i\pi/4}(X-1)-4+\pb\left(\left|-4\sqrt X -2\sqrt\pi e^{-3i\pi/4}(X-1)+4\right|+4/3\right).
\end{align*}
The part inside of the $\pb$ symbol is shown to be a decreasing function of $X$ over $[0,1]$, and hence its maximum value over $[0,1]$ is equal to $|2\sqrt\pi e^{-3i\pi/4}+4|+4/3\leq 16/3$, which concludes.

We now turn to the last estimate of the lemma. We have 
\[\Delta_3(X)=\int_0^{X}(X-x)\left(-2x^{-1/2}+\int_0^x \frac{e^{-it}-1}{t^{3/2}}\diff t\right)\diff x=\int_0^X\Delta_2(x)\diff x.\]
By the computations above, this concludes.
\end{proo}
\begin{lemm}{}\label{lemmePgoth}
Let $d$ be an integer, $\varepsilon=c/\log T$ with $c>0$, and let $0\leq v\leq H$. We define 
\[\mathfrak P_d(s)=\prod_{p|d}(1-p^{-(1+s+iv)})^{-1/2}.\]
Put, for $r=2$, $3$, 
\[\mathfrak P_d(s)=\sum_{j=0}^{r-1}\frac{s^j}{j!} \mathfrak P_d^{(j)}(0)+s^r R_{r,d}(s).\]
Then for $-1\leq t\leq 1$, $j=0$, $1$, $2$ we have:
\[\mathfrak P_d^{(j)}(0)=\O\left(\prod_{p|d}(1+p^{-3/4})\right),\qquad R_{r,d}(\varepsilon+it)=\O\left(\prod_{p|d}(1+p^{-3/4})\right).\]
\end{lemm}
\begin{proo}{}
We consider the rectangle $Rec$ defined by the lines $\Re(s)=\varepsilon\pm 1/5$ and $\Im(s)=\pm 1$ where $\mathfrak P_d(s)$ is bounded by $\prod_{p|d}(1-p^{-1+1/5})^{-1/2}\ll \prod_{p|d}(1+p^{-3/4})$. By Cauchy's integral formula, we find 
\[\mathfrak P_d^{(j)}(0)=\frac{j!}{2\pi i}\int_{Rec}\frac{\mathfrak P_d(z)}{z^{j+1}}\diff z=\O\left(\prod_{p|d}(1+p^{-3/4})\right)\]
and
\[R_{r,d}(\varepsilon +it)=\frac{1}{2\pi i}\int_{Rec}\frac{\mathfrak P_d(z)}{(z-(\varepsilon+it))z^r}\diff z=\O\left(\prod_{p|d}(1+p^{-3/4})\right).\]
\end{proo}

\begin{lemm}{}\label{lemm12}
Let $d,\delta,q\geq 1$ be a positive integers, $(d,q)=1$, and $0< v\leq 1/\sqrt{\log \xi}$. Then for $r=2$, $3$, we have 
\begin{align*}
\sum_{\substack{(n,qd)=1\\n\delta<\xi}}\frac{\tau_{-1/2}(n)\log^{r-1}\left(\frac{\xi}{n\delta}\right)}{n^{1+iv}}&=-\frac{(r-1)!}{2\sqrt\pi}\prod_{p|qd}(1-p^{-(1+iv)})^{-1/2}v^{3/2-r}\Delta_{r}(\log(\xi/\delta)v)
\\&\qquad +\O\left(\prod_{p|d}(1+p^{-3/4})\left((\log T)^{r-5/2}+(\log\log T)^{r-3/2)}\right)\right).
\end{align*}
\end{lemm}
\begin{proo}{}
The lemma being true if $\xi/\delta<2$, we may suppose that $\xi/\delta\geq 2$. Put $\varepsilon=c/\log T$ with $c>0$. For $x>0$ and $r=2$, $3$, we have the Mellin transform
\[\frac{1}{2\pi i}\int_{\varepsilon-i\infty}^{\varepsilon+i\infty}\frac{x^s}{s^r}\diff s=\left\{\begin{array}{cl}
0&\text{if }0<x\leq 1,
\\\frac{(\log x)^{r-1}}{(r-1)!}& \text{if }x>1.
\end{array}\right.\]
Therefore, as in \nameref{lemmeA},
\begin{align*}
\sum_{\substack{(n,qd)=1\\ n\delta<\xi}}&\frac{\tau_{-1/2}(n)\log^{r-1}\left(\frac{\xi}{n\delta}\right)}{n^{1+iv}}=\frac{(r-1)!}{2\pi i}\int_{\varepsilon-i\infty}^{\varepsilon+i\infty}\frac{\left(\frac{\xi}\delta\right)^s}{s^r}\sum_{(n,dq)=1}\frac{\tau_{-1/2}(n)}{n^{1+s+iv}}\diff s
\\&=\frac{(r-1)!}{2\pi i}\int_{\varepsilon-i/\log\log T}^{\varepsilon+i/\log\log T}\frac{\left(\frac{\xi}\delta\right)^s}{s^r}Z(1+s+iv)\diff s+\O\left(\sqrt{\frac{d}{\phi(d)}}(\log\log T)^{r-3/2}\right)
\\&=\numberthis\label{eqfrakPint}\frac{(r-1)!}{2\pi i}\int_{\varepsilon-i}^{\varepsilon+i}\frac{\left(\frac{\xi}\delta\right)^s}{s^r}\mathfrak P_{qd}(s)\sqrt{s+iv}\diff s+\O\left(\sqrt{\frac{d}{\phi(d)}}\left((\log T)^{r-5/2}+(\log\log T)^{r-3/2)}\right)\right).
\end{align*}
Note that, by Lemma 9 of \cite{selberg_zeta}, we have for, $k=1,2,3$:
\[\frac{1}{2\pi i}\int_{\varepsilon-i}^{\varepsilon+i}\left(\frac{\xi}\delta\right)^s\frac{\sqrt{s+iv}}{s^{k}}\diff s=\O\left(\log(T)^{k-3/2}+1\right).\]
We may now apply this fact, \nameref{intracine} and \nameref{lemmePgoth} to get
\begin{align*}
&\frac{1}{2\pi i}\int_{\varepsilon-i}^{\varepsilon+i}\frac{\left(\frac{\xi}\delta\right)^s}{s^r}\mathfrak P_{qd}(s)\sqrt{s+iv}\diff s=\sum_{j=0}^{r-1}\frac{1}{2\pi i}\frac{\mathfrak P_{qd}^{(j)}(0)}{j!}\int_{\varepsilon-i}^{\varepsilon+i}\left(\frac{\xi}\delta\right)^s\frac{\sqrt{s+iv}}{s^{r-j}}\diff s
\\&\qquad+\frac{1}{2\pi i}\int_{\varepsilon-i}^{\varepsilon+i}\left(\frac{\xi}\delta\right)^s\sqrt{s+iv}R_{r,qd}(s)\diff s
\\&=-\mathfrak P_{qd}(0)\frac{v^{3/2-r}}{2\sqrt \pi}\Delta_{r}(\log(\xi/\delta)v)+\sum_{j=1}^{r-1}\mathfrak P_{qd}^{(j)}(0)\O\left((\log T)^{r-j-3/2}+1\right)
\\&\qquad +\O\left(\prod_{p|d}(1+p^{-3/4})\right)
\\&\numberthis\label{eqfrakP}=\mathfrak -\mathfrak P_{qd}(0)\frac{v^{3/2-r}}{2\sqrt \pi}\Delta_{r}(\log(\xi/\delta)v)+\O\left(\prod_{p|d}(1+p^{-3/4})\left(\log(T)^{r-5/2}+1\right)\right).
\end{align*}
Since $\sqrt{d/\phi(d)}\ll \prod_{p|d}(1+p^{-3/4})$, one may use \eqref{eqfrakPint} and \eqref{eqfrakP} to conclude.
\end{proo}
\begin{coro}{}\label{coro-molli2}
Let $d,q\geq 1$ be co-prime positive integers, $1\leq \delta\leq \xi$ be such that $\delta|d^\infty$, and $0<v\leq 1/\log \xi$. Then
\begin{align*}
\sum_{(n,qd)=1}\frac{\M(n\delta)\tau_{-1/2}(n\delta)}{n}\frac{\sin(v\log(n\delta))}{v}&=\pb\Biggl[\frak C_4\tau_{1/2}(\delta)\sqrt{\frac{dq}{\phi(d)\phi(q)}}(\log \xi)^{1/2}\Biggr]
\\&\qquad+\O\left(\tau_{1/2}(\delta)\prod_{p|d}(1+p^{-3/4})(\log T)^{1/3}\right),
\end{align*}
where 
\begin{align*}
&\frak C_4:=
\\&\frac{1}{3(1-\theta)\sqrt\pi}\left[\frac{9(1+\theta^{5/2})}{5}+2\max\{3\sqrt{1-\theta},\min\{ \max\{1,|1-4\theta|\}+3\theta\sqrt\theta, 1-\theta+3\theta\sqrt{1-\theta}\}\}\right] .
\end{align*}
\end{coro}
\begin{proo}{}
Multiplying the formula of \nameref{lemm12} for $r=2$ by $\log \xi$ and subtracting the formula for $r=3$, we find that for any $0<u\leq v$:
\begin{align*}
&\sum_{\substack{(n,qd)=1\\ n\delta<\xi}}\frac{\tau_{-1/2}(n)\log\left(\frac{\xi}{n\delta}\right)}{n^{1+iu}}\log(n\delta)=\frac{-1}{2\sqrt\pi}\prod_{p|qd}(1-p^{-(1+iu)})^{-1/2}\Biggl((\log\xi) u^{-1/2}\Delta_2(\log(\xi/\delta)u)-
\\& 2u^{-3/2}\Delta_{3}(\log(\xi/\delta)u)\Biggr)+\O\left(\prod_{p|d}(1+p^{-3/4})\log T\sqrt{\log\log T}\right).
\end{align*}
Thus by \nameref{lemm_esti_delta},
\begin{align*}
&\sum_{\substack{(n,qd)=1\\ n\delta<\xi}}\frac{\tau_{-1/2}(n)\log\left(\frac{\xi}{n\delta}\right)}{n^{1+iu}}\log(n\delta)=\O\left(\prod_{p|d}(1+p^{-3/4})(\log T)^{4/3}\right)
\\&+\pb\Biggl[\frac{18}{15\sqrt\pi}\sqrt{\frac{dq}{\phi(d)\phi(q)}} u(\log \xi)^{5/2}\Biggr]-\frac{2}{3\sqrt\pi}\prod_{p|qd}(1-p^{-(1+iu)})^{-1/2}(\log\xi-4\log \delta)\sqrt{\log(\xi/\delta)}.
\end{align*}
Using the relation between $\M_{sel}$ and $\M$ (see \eqref{lien-molli}), we find that 
\begin{align*}
\numberthis\label{eq_a_mult}&\sum_{(n,qd)=1}\frac{\M(n\delta)\tau_{-1/2}(n\delta)}{n^{1+iu}}\log(n\delta)=\pb\Biggl[\frac{18(1+\theta^{5/2})}{(1-\theta)15\sqrt\pi}\tau_{1/2}(\delta)\sqrt{\frac{dq}{\phi(d)\phi(q)}} u(\log \xi)^{3/2}\Biggr]
\\&+\pb\left[\frac{2\tau_{1/2}(\delta)}{(1-\theta)3\sqrt\pi}\sqrt{\frac{dq}{\phi(d)\phi(q)}}\left|\left(1-4\frac{\log \delta}{\log \xi}\right)\sqrt{\log\frac{\xi}\delta}-\left(\theta-4\frac{\log \delta}{\log \xi}\right)\sqrt{\log^+\left(\frac{\xi^\theta}{\delta}\right)}\right|\right]
\\&+\O\left(\tau_{1/2}(\delta)\prod_{p|d}(1+p^{-3/4})(\log T)^{1/3}\right).
\end{align*}
We denote by $\mathcal Q$ the quantity in absolute value, in the second $\pb$. Then, if $\delta>\xi^\theta$, observe that for any $\theta\leq x\leq 1$, we have $|1-4x| \leq 3$, and hence
\begin{align*}
|\mathcal Q|\leq 3\sqrt{\log\frac{\xi}{\delta}}\leq 3\sqrt{1-\theta}\sqrt{\log \xi}.
\end{align*}
If $\delta\leq \xi^\theta$, then observe that for $x\in[0,\theta]$, we have $|1-4x|\leq \max(1,|1-4\theta|)$ and $|\theta-4x|\leq 3\theta$, and hence
\begin{align*}
|\mathcal Q|\leq \left(\max(1,|1-4\theta|)+ 3\theta	\sqrt\theta\right)\sqrt{\log\xi}.
\end{align*}
We may also write the following, since $\sqrt a-\sqrt b\leq \sqrt{a-b}$ for any real numbers $0\leq a<b$:
\begin{align*}
|\mathcal Q|&\leq \left|1-4\frac{\log \delta}{\log \xi}-\theta+4\frac{\log \delta}{\log \xi}\right|\sqrt{\log\frac{\xi}\delta}+\left|\theta-4\frac{\log \delta}{\log \xi}\right|\left|\sqrt{\log \frac{\xi^\theta}{\delta}}-\sqrt{\log\frac{ \xi}{\delta}}\right|
\\&\leq (1-\theta)\sqrt{\log \xi}+3\theta\sqrt{1-\theta}\sqrt{\log \xi}.
\end{align*}
Therefore we always have
\[|\mathcal Q|\leq \max\{3\sqrt{1-\theta},\min\{ \max\{1,|1-4\theta|\}+3\theta\sqrt\theta, 1-\theta+3\theta\sqrt{1-\theta}\}\sqrt{\log \xi}\}=:\C(\theta)\sqrt{\log \xi}.\]

Multiplying \eqref{eq_a_mult} by $\delta^{-iu}$, integrating with respect to $u$ from $0$ to $v$, and dividing by $v$, we find that
\begin{align*}
&\sum_{(n,qd)=1}\frac{\M(n\delta)\tau_{-1/2}(n\delta)}{nv}i\left((n\delta)^{-iv}-1\right)=\pb\Biggl[\frac{9(1+\theta^{5/2})}{(1-\theta)15\sqrt\pi}\tau_{1/2}(\delta)\sqrt{\frac{dq}{\phi(d)\phi(q)}} v(\log \xi)^{3/2}\Biggr]
\\&+\pb\left[\frac{2\tau_{1/2}(\delta)}{(1-\theta)3\sqrt\pi}\sqrt{\frac{dq}{\phi(d)\phi(q)}}\C(\theta)\sqrt{\log\xi}\right]+\O\left(\tau_{1/2}(\delta)\prod_{p|d}(1+p^{-3/4})(\log T)^{1/3}\right).
\end{align*}
Now by taking the real part of this we get the required result, since $v\log\xi\leq 1$.
\end{proo}

\section{Interlude: integral mean value theorems for Dirichlet $L$-functions}\label{sec:inter}
Before turning to the proofs of \eqref{estiI} and \nameref{lemme-mollifer} -which are quite long and much harder than the previous one-, we need to state some results about integral moments of $L$-functions. We will be following Bauer's papers \cite{bauer2},\cite{bauer}, where he essentially proved our estimates (he studied a slightly different integral, and only in the case $\alpha\neq \beta$ in the following theorem). Therefore we will detail how to get the main term in the following theorem -whether $\alpha=\beta$ or not-, but we will mostly rely on Bauer's papers to estimate the error terms. 

Note that this section differs from Selberg's original proof, and this allows us to simplify the proof a little and get stronger errors terms in crucial estimates. The consequence of that is that we may take our mollifier to be a bit longer than the one used by Selberg: the length of his mollifier was about $T^{1/10}$, while ours has size of roughly $T^{1/8}$.
\begin{theo}{see Section 3 of \cite{bauer}}\label{th:meanvalue} We put $\mathcal T:=qT/(2\pi)$. Let $q\geq 1$, $\chi$ be an even primitive Dirichlet character mod $q$. Let $\alpha=iU/\log T,\beta =iV/\log T$ be distinct complex numbers, where $U$ and $V$ are real numbers bounded by an absolute constant. Let $h,k\leq T^\omega$, $0<\omega<1/2$, be positive integers and fix $1/2\leq c<1$. We have 
\begin{align*}
\frac1{i}&\int_{c+i}^{c+iT}L(s+\alpha,\chi)L(1-s-\beta,\overline \chi)h^{-s}k^{s-1}\diff s=\chi\left(\frac{k}{(h,k)}\right)\overline\chi\left(\frac{h}{(h,k)}\right)T
\\&\times \Biggl(L(1-\alpha+\beta,\chi_0)\frac{\mathcal T^{\beta-\alpha}(h,k)^{1-\alpha+\beta}}{(1-\alpha)h^{1-\alpha}k^{1+\beta}}+ L(1+\alpha-\beta,\chi_0)\frac{(h,k)^{1+\alpha-\beta}}{(1-\beta)h^{1-\beta}k^{1+\alpha}}\Biggr)+R_c(T;h,k;\alpha,\beta),
\end{align*}
where $R_c(T;h,k;\alpha,\beta)$ is an error term such that for any small, but fixed, $\varepsilon>0$ and any bounded complex sequence $(\theta(n))$:
\[\numberthis\label{condition}\sum_{h,k\leq T^\omega}\theta(h)\overline{\theta(k)}R_c(T;h,k;\alpha,\beta)\ll (T^{1/2+\omega}+T^{1/3+\omega(2-c)})T^\varepsilon.\]
If we suppose that $\alpha=\beta$, then result holds if the main term is replaced by 
\[\chi\left(\frac{k}{(h,k)}\right)\overline\chi\left(\frac{h}{(h,k)}\right)\frac{(h,k)}{h^{1-\alpha}k^{1+\alpha}}\frac{T}{1-\alpha}\frac{\phi(q)}{q}\left(\log\left(\frac{T (h,k)^2}{hk}\right)+\frak Q(q,\alpha)\right),\]
for a quantity $\frak Q(q,\alpha)$ defined in \eqref{defQ} below.
\end{theo}
\begin{proo}{}
Let $M:=\frac1{i}\int_{c+i}^{c+iT}L(s+\alpha,\chi)L(1-s-\beta,\overline \chi)h^{-s}k^{s-1}\diff s$. We move the line of integration to the right by the residue theorem, say at $c'=1+\varepsilon$ for some small $\varepsilon>0$. Using the fact that $L(1/2+\sigma+it,\chi)\ll t^{(1-\sigma)/3+o(1)}$ if $t>1$ and $\sigma\in[1/2,1)$ (see Theorem 2 of \cite{kol}), we get
\[\numberthis\label{eqM1}M=\frac 1i\int_{c'+i}^{c'+iT}L(s+\alpha,\chi)L(1-s-\beta,\overline \chi)h^{-s}k^{s-1}\diff s+\O(T^{1/3+\varepsilon} k^\varepsilon/h^c).\]
We define 
\[\numberthis\label{eq_err0} Err_0:=T^{1/3+\varepsilon} k^\varepsilon/h^c,\]
and we write the functional equation in the form 
\[\numberthis\label{eq:funceq} L(s,\chi)=\Xi(s,\chi)L(1-s,\overline \chi)\]
We will use Lemma 1 of \cite{bauer}, which states the following: Let $r>0$, and for $x\in\RR$ let $e(x)=\exp(2\pi ix)$. We also introduce 
\[E_{c'}(r,T)=T^{c'-1/2}+\frac{T^{c'+1/2}}{|T-r|+T^{1/2}},\qquad G(\chi)=\sum_{n=1}^q \chi(n) e\left(\frac{n}q\right).\]
Then this lemma reveals that if $r\leq \mathcal T$, then 
\begin{align*}
\int_{c'+i}^{c'+iT}\Xi(1-s-\beta,\overline{\chi})r^{-s}\diff s&=2\pi iG(\overline \chi)\mathcal T^\beta q^{-1}e\left(-\frac rq\right)
+\O\left(r^{-1}+r^{-c'}E_{c'}\left(\frac{2\pi r}q,T\right)\right),
\end{align*}
and if $r>\mathcal T$, then
\[\int_{c'+i}^{c'+iT}\Xi(1-s-\beta,\overline{\chi})r^{-s}\diff s=\O\left(r^{-c'}E_{c'}\left(\frac{2\pi r}q,T\right)\right).\]
Using this lemma together with the functional equation \eqref{eq:funceq} in \eqref{eqM1}, we may write that 
\begin{align*}
M&=\frac{1}i\int_{c'+i}^{c'+iT}L(s+\alpha,\chi)L(s+\beta,\chi)\Xi(1-s-\beta,\overline\chi)h^{-s}k^{s-1}\diff s+\O(Err_0)
\\&=\frac{1}{ik}\sum_{n,m\geq 1}\frac{\chi(n)\chi(m)}{n^\alpha m^\beta}\int_{c'+i}^{c'+iT}\left(\frac{k}{hnm}\right)^s\Xi(1-s-\beta,\overline\chi)\diff s+\O(Err_0)
\\\numberthis\label{expM}&=\frac{2\pi G(\overline \chi)\mathcal T^\beta}{qk}\sum_{nm\leq \mathcal T k/h}\frac{\chi(n)\chi(m)}{n^\alpha m^\beta}e\left(-\frac{nmh}{qk}\right)+\O\left(Err_0+Err_1(h,k)\right)
\end{align*}
where 
\[\numberthis\label{eq_err1} Err_1(h,k)=\sum_{m,n\geq 1}\frac1{n^{\Re(\alpha)}m^{\Re(\beta)}}\frac{k^{\varepsilon}}{(nmh)^{1+\varepsilon}}E_{c'}(2\pi nmh/(qk),T)+\sum_{m,n:mnh/k\leq \mathcal T}\frac{1}{mnh}.\]

Now we use Perron's formula (see (2.7) of \cite{bauer2}) to write that for any large $x\geq 1$ and for $\varrho=\max(|\alpha|,|\beta|)$:
\begin{align*}
\numberthis\label{perron}\sum_{nm\leq \mathcal T k/h}&\frac{\chi(n)\chi(m)}{n^\alpha m^\beta}e\left(-\frac{nmh}{qk}\right)=\frac{1}{2\pi i}\int_{c'-ix}^{c'+ix}\frac{(\mathcal T k/h)^s}s\sum_{n,m}\frac{\chi(n)\chi(m)}{n^{s+\alpha}m^{s+\beta}}e\left(-\frac{nmh}{qk}\right)\diff s
\\&+\O\left(\frac{(\mathcal T k/h)^{1+\varepsilon}}{x}+\frac{2^{2\varrho}(\mathcal T k/h)^{1+2\varrho}\log (2\mathcal T k/h)}{x}+(2\mathcal T k/h)^{2\varrho}\right).
\end{align*}
This error term is easily shown to be 
\[\numberthis\label{eq_err2}\ll\frac{(Tk/h)^{1+\varepsilon}}x+1=:Err_2(x,h,k).\]
Now we denote the Hurwitz zeta-function by $\zeta(s,y)$, $0<y\leq 1$, and we let $H':=h/(h,k)$ and $K:=k/(h,k)$. The Hurwitz zeta-function $\zeta(s,y)$ is the analytic continuation of the function $\sum_{n=0}^\infty (n+y)^{-s}$, defined for $\Re(s)>1$, to the whole complex plan except at $s=1$. When $s$ is close to $1$, we have $\zeta(s,y)=1/(s-1)+\O_y(1)$, and hence $\zeta(s,y)-\zeta(s)$ is regular on the whole complex plan. We define 
\[D(s):=\sum_{\nu,\mu\leq qK}\chi(\nu)\chi(\mu)e\left(-\frac{\nu \mu h}{qk}\right)\zeta\left(s+\alpha,\frac{\nu}{qK}\right)\zeta\left(s+\beta,\frac{\mu}{qK}\right),\]
and
\[D^*(s):=\sum_{\nu,\mu\leq qK}\chi(\nu)\chi(\mu)e\left(-\frac{\nu \mu H'}{qK}\right)\left(\zeta\left(s+\alpha,\frac{\nu}{qK}\right)-\zeta(s+\alpha)\right)\left(\zeta\left(s+\beta,\frac{\mu}{qK}\right)-\zeta(s+\beta)\right).\]
For $s\in\CC$, we have 
\begin{align*}
\sum_{n,m}\frac{\chi(n)\chi(m)}{n^{s+\alpha}m^{s+\beta}}e\left(-\frac{nmh}{qk}\right)&=\sum_{1\leq \nu,\mu\leq qK}\chi(\nu)\chi(\mu)e\left(-\frac{\nu\mu h}{qk}\right)\sum_{\substack{n\equiv \nu\mod{qK}\\m\equiv \mu\mod{qK}}}n^{-s-\alpha}m^{-s-\beta}
\\\numberthis\label{eqDs}&=(qK)^{-2s-\alpha-\beta}D(s).
\end{align*}
Bauer shows in pages 29-30 of \cite{bauer2} that $D^*(s)=D(s)-E_1(s)-E_2(s)+E_3(s)$, with 
\[E_1(s)=Kq^{s+\beta}G(\chi)\chi(K)\overline\chi(-H')L(s+\beta,\chi_0)\zeta(s+\alpha),\]
\[E_2(s)=Kq^{s+\alpha}G(\chi)\chi(K)\overline{\chi}(-H')L(s+\alpha,\chi_0)\zeta(s+\beta),\]
\[E_3(s)=K\phi(q)G(\chi)\chi(K)\overline\chi(-H')\zeta(s+\alpha)\zeta(s+\beta),\]
where $\phi$ is Euler's totient function and $\chi_0$ is the principal character modulo $q$. Using the fact that $\chi$ is assumed to be even, we find that $\overline\chi(-H')=\overline\chi(H')$. Let $\R_{c'}(x)$ be the rectangle with vertices at $c'\pm ix$ and $1/2\pm ix$, and define $\Gamma_{c'}(x)$ as the path along the upper, left and lower part of $\R_{c'}(x)$. Then by using \eqref{eqDs} the main term of \eqref{perron} becomes 
\[\numberthis\label{eq_res_rec}\frac{1}{2\pi i}\int_{\R_{c'(x)}}\frac{(\mathcal T k/h)^s}s\frac{D^*(s)+E_1(s)+E_2(s)-E_3(s)}{(qK)^{2s+\alpha+\beta}}\diff s-\frak I(h,k),\]
with
\[\frak I(h,k):=\frac{1}{2\pi i}\int_{\Gamma_{c'(x)}}\frac{(\mathcal T k/h)^s}s\sum_{n,m}\frac{\chi(n)\chi(m)}{n^{s+\alpha}m^{s+\beta}}e\left(-\frac{nmh}{qk}\right)\diff s.\]

First, we suppose that $\alpha\neq \beta$. Since $\zeta(s,y)-\zeta(s)$ is regular for $0<y\leq 1$ -and hence $D^*$ is regular-, and since $\alpha\neq \beta$, the only poles of the first integral are simple poles at $1-\alpha$ and $1-\beta$. To compute the residues involved, note that $L(s,\chi_0)=\zeta(s)\prod_{p|q}(1-p^{-s})$ and hence $\Res_{s= 1}(L(s,\chi_0))=\phi(q)/q$. Thus, $\Res_{s=1-\alpha}(E_2(s))=\Res_{s=1-\alpha}(E_3(s))$ which implies that 
\[\Res_{s=1-\alpha}(E_1(s)+E_2(s)-E_3(s))=\Res_{s=1-\alpha}(E_1(s))=K q^{1-\alpha+\beta}G(\chi)\chi(K)\overline\chi(H')L(1-\alpha+\beta,\chi_0).\]
The pole at $s=1-\beta$ is handled similarly, and hence the Residue Theorem yields that the \eqref{eq_res_rec} integral along the rectangle is equal to 
\begin{align*}
\frac{G(\chi)\chi(K)\overline\chi(H')}{qK}\frac{\mathcal T k}h\Biggl( L(1-\alpha+\beta,\chi_0)\frac{(\mathcal T k/h)^{-\alpha}K^{\alpha-\beta}}{1-\alpha}+ L(1+\alpha-\beta,\chi_0)\frac{(\mathcal T k/h)^{-\beta}K^{-\alpha+\beta}}{1-\beta}\Biggr).
\end{align*}

Putting these back in \eqref{perron} and \eqref{expM}, and using the fact that $G(\chi)G(\overline \chi)=q$, we get
\begin{align*}
\numberthis\label{finM}M&=\chi(K)\overline\chi(H')\frac T{hK}\left(L(1-\alpha+\beta,\chi_0)\frac{\mathcal T^{\beta-\alpha}(k/h)^{-\alpha}K^{\alpha-\beta}}{1-\alpha}+ L(1+\alpha-\beta,\chi_0)\frac{(k/h)^{-\beta}K^{-\alpha+\beta}}{1-\beta}\right)
\\&-\frac{2\pi}{qk}G(\overline{\chi})\mathcal T^\beta\frak I(h,k)+\O\left(Err_0+Err_1(h,k)+\frac{Err_2(x,h,k)}k\right).
\end{align*} 

Now, we suppose that $\alpha=\beta$. We see that the only difference lies in the fact that the only pole of the integral \eqref{eq_res_rec} is now of order $2$, at $s=1-\alpha$. Let us compute the corresponding residue, that we denote by $R$. We have
\begin{align*}
R&=\frac{\diff}{\diff s}\Biggl[\frac{E_1(s)+E_2(s)-E_3(s)}{s(qK)^{2s+2\alpha}}\left(\frac{\mathcal T k}{h}\right)^s(s-1+\alpha)^2\Biggr]_{s=1-\alpha}
\\&=\numberthis\label{eq:residue} \frac{G(\chi)\chi(K)\overline\chi(H')}{q^2K} \left(\frac{\mathcal T k}{h}\right)^{1-\alpha}\frac{\diff}{\diff z}\Biggl[ \frac{2q^z \prod_{p|q}(1-p^{-z})-\phi(q)}{(z-\alpha)(qK)^{2(z-1)}}\zeta^2(z)\left(\frac{\mathcal T k}{h}\right)^{z-1}(z-1)^2\Biggr]_{z=1}.
\end{align*} 
Recall that $\mathcal T=qT/(2/\pi)$. The derivative part of this is equal to 
\begin{align*}
&\frac{\left(2\phi(q)\log q +2q\sum_{p_1|q}\prod_{\substack{p_2|q\\p_1\neq p_2}}(1-p_2^{-1})\frac{\log p_1}{p_1}\right)(1-\alpha)-\phi(q)\left(1+2(1-\alpha)\log(qK)\right)}{(1-\alpha)^2}
\\&\qquad+\frac{\phi(q)}{(1-\alpha)}\left[\log\left(\frac{\mathcal T k}h\right)+\Res_{s=1}(\zeta(s)^2(s-1)^2)\right]
\\&=\frac{\phi(q)}{1-\alpha}\left(2\log q+2\sum_{p|q}\frac{\log p}{p-1}-\frac{1}{1-\alpha}-2\log(qK)+\log\left(\frac{\mathcal T k}h\right)+2\gamma\right).
\\&=\frac{\phi(q)}{1-\alpha}\left(\log\left(\frac{T k}{hK^2}\right)+\frak Q(q,\alpha)\right),
\end{align*}
where $\gamma$ is the usual Euler-Mascheroni constant and
\[\numberthis\label{defQ}\frak Q(q,\alpha):=-\frac{1}{1-\alpha}+\log\frac{q}{2\pi}+2\gamma+2\sum_{p|q}\frac{\log p}{p-1}.\]

Putting this in \eqref{eq:residue}, the Residue Theorem yields that the \eqref{eq_res_rec} integral along the rectangle is equal to 
\[G(\chi)\chi(K)\overline\chi(H')\left(\mathcal T\frac k h\right)^{-\alpha}\frac{k}{qhK}\frac{T}{2\pi}\frac{\phi(q)}{1-\alpha}\left(\log\left(\frac{T (h,k)^2}{hk}\right)+\frak Q(q,\alpha)\right),\]
Again, using the equations \eqref{perron} and \eqref{expM} above, we have proved that
\begin{align*}
M&=\chi(K)\overline\chi(H')\frac{(h,k)}{h^{1-\alpha}k^{1+\alpha}}\frac{T}{1-\alpha}\frac{\phi(q)}{q}\left(\log\left(\frac{T (h,k)^2}{hk}\right)+\frak Q(q,\alpha)\right)-\frac{2\pi}{qk}G(\overline{\chi})\mathcal T^\alpha\frak I(h,k)
\\&+\O\left(Err_0+Err_1(h,k)+\frac{Err_2(x,h,k)}k\right).
\end{align*} 
One should note that we find the same main term found by Matsumoto \cite{matsumoto}. 

Now we turn the error terms. We have by definition (see\eqref{eq_err0} and \eqref{eq_err2})
\[\sum_{h,k\leq T^\omega}|Err_0|\ll T^{1/3+\omega(2-c)+2\varepsilon},\qquad \sum_{h,k\leq T^\omega} \frac{|Err_2(x,h,k)|}k\ll T^{2\omega}+\frac{T^{1+\omega+\varepsilon}}{x}.\]
By Lemma 2.2.3 of \cite{bauer2}, the definition \eqref{eq_err1} of $Err_1$ implies that
\[\sum_{h,k\leq T^\omega}|Err_1(h,k)|\ll T^{1/2+\omega+2\varepsilon}.\]

Now it only remains to estimate the contribution of $\frak I(h,k)/k$ in \eqref{finM}. If $x$ is chosen as a large power of $T$, then Section 5 of \cite{bauer} shows that for any bounded complex sequence $(\theta(n))$,
\[\sum_{h,k\leq T^\omega}\theta(h)\overline{\theta(k)}\frac{\frak I(h,k)}k\ll T^{1/2+\omega+2\varepsilon}.\]
\end{proo}

From this, we easily deduce three corollaries. We keep the notation for $R_c(T;h,k;\alpha,\beta)$. The first corollary is straightforward by taking $\alpha=\beta=0$:
\begin{coro}{}\label{coroegal}
Let $h,k\leq T^\omega$, $0<\omega<1/2$ be integers. Then
\begin{align*}
\frac{1}{\sqrt{hk}}\int_T^{2T} |L(1/2+it,\chi)|^2\left(\frac hk\right)^{-it}\diff t&=\chi\left(\frac{k}{(h,k)}\right)\overline\chi\left(\frac{h}{(h,k)}\right)(h,k)\frac{\phi(q)}{q}\frac{T}{kh}
\\&\times\left(\log\left(\frac{T(h,k)^2}{hk}\right)+\frak Q(q,0)+2\log 2\right)+R_{1/2}(T;h,k;0,0).
\end{align*}
\end{coro}
The second one is proved by partial summation.
\begin{coro}{}\label{coroipp}
Let $\alpha=iU/\log T,\beta= iV/\log T$ be distinct complex numbers, where $U$ and $V$ are real numbers bounded by an absolute constant. Let $h,k\leq T^\omega$, $0<\omega<1/2$, be positive integers. For any smooth function $g:(0,\infty)\to\CC$, we have:
\begin{align*}
&\qquad\frac{1}{\sqrt{hk}}\int_{1}^{T}L(1/2+it+\alpha,\chi)L(1/2-it-\beta,\overline \chi)\left(\frac{h}k\right)^{-it}g(t)\diff t=\frac{\chi\left(\frac{k}{(h,k)}\right)\overline\chi\left(\frac{h}{(h,k)}\right)}{hk}
\\&\times \int_{1}^{T}g(t)\Biggl(L(1-\alpha+\beta,\chi_0)\frac{(1+\beta-\alpha)(qt/2\pi)^{\beta-\alpha}(h,k)^{1-\alpha+\beta}}{(1-\alpha)h^{-\alpha}k^{\beta}}+ L(1+\alpha-\beta,\chi_0)\frac{(h,k)^{1+\alpha-\beta}}{(1-\beta)h^{-\beta}k^{\alpha}}\Biggr)\diff t
\\&\qquad+R(T;h,k;\alpha,\beta;g),
\end{align*} 
where $R(T;h,k;\alpha,\beta;g)$ is an error term satisfying
\[\sum_{h,k\leq T^\omega}\theta(h)\overline{\theta(k)}R(T;h,k;\alpha,\beta;g)\ll (T^{1/2+\omega}+T^{1/3+3\omega/2})T^\varepsilon\left(|g(T)|+\int_1^T|g'(t)|\diff t+1\right)\]
for any $\varepsilon>0$ and any bounded complex sequence $(\theta(n))$.
\end{coro}
\begin{proo}{}
\nameref{th:meanvalue} proves that for all $T>1$ and any distinct complex numbers $\alpha,\beta$, we have
\begin{align*}
&F(T):=\frac{1}{\sqrt{hk}}\int_1^TL(1/2+it+\alpha,\chi)L(1/2-it-\beta,\overline\chi)\left(\frac hk\right)^{-it}\diff t=\frac{\chi\left(\frac{k}{(h,k)}\right)\overline\chi\left(\frac{h}{(h,k)}\right)}{hk}T
\\&\numberthis\label{eq:ipp}\times\left(L(1-\alpha+\beta,\chi_0)\frac{(qT/2\pi)^{\beta-\alpha}(h,k)^{1-\alpha+\beta}}{(1-\alpha)h^{-\alpha}k^\beta}+L(1+\alpha-\beta,\chi_0)\frac{(h,k)^{1+\alpha-\beta}}{(1-\beta)h^{-\beta}k^\alpha}\right)+R(T),
\end{align*}
where $R(T)$ is an error term which is $R_{1/2}(T;h,k;\alpha,\beta)$ when $\alpha$ and $\beta$ are of the form $iU/\log T$ and $iV/\log T$. Thus if $g$ is a smooth function over $(0,\infty)$, we get by partial integration
\begin{align*}
&\frac{1}{\sqrt{hk}}\int_1^TL(1/2+it+\alpha,\chi)L(1/2-it-\beta,\overline\chi)\left(\frac hk\right)^{-it}g(t)\diff t=g(T)F(T)-\int_1^Tg'(t)F(t)\diff t
\\&=\frac{\chi\left(k/(h,k)\right)\overline\chi\left(h/(h,k)\right)}{hk}\int_1^T g(t)\Biggl(L(1-\alpha+\beta,\chi_0)(1+\beta-\alpha)\frac{(qt/2\pi)^{\beta-\alpha}(h,k)^{1-\alpha+\beta}}{(1-\alpha)h^{-\alpha}k^\beta}
\\&\qquad+L(1+\alpha-\beta,\chi_0)\frac{(h,k)^{1+\alpha-\beta}}{(1-\beta)h^{-\beta}k^\alpha}\Biggr)\diff t+g(T)R(T)-\int_1^Tg'(t)R(t)\diff t+\O(1).
\end{align*}
Now, we take our $\alpha$ and $\beta$ to be of the wanted form, we let $(\theta(n))$ be a bounded sequence of complex numbers and we fix $\varepsilon>0$. Then by property of $R(T;h,k;\alpha,\beta)=R(T)$:
\begin{align*}
&\sum_{h,k\leq T^\omega}\theta(h)\overline{\theta(k)}\left[g(T)R(T)-\int_1^Tg'(t)R(t)\diff t+\O(1)\right]
\\&\ll (T^{1/2+\omega}+T^{1/3+3\omega/2})T^\varepsilon |g(T)|+\int_1^T|g'(t)|(t^{1/2+\omega}+t^{1/3+3\omega/2})t^\varepsilon\diff t+T^{2\omega}
\\&\ll (T^{1/2+\omega}+T^{1/3+3\omega/2})T^\varepsilon\left(|g(T)|+\int_1^T|g'(t)|\diff t+1\right).
\end{align*}
This concludes.
\end{proo}
In the case where $g(t)=(2\pi/t)^{i v/2}$ for some $v=V/\log T$, a choice motivated by Stirling's formula as we shall see later, we find the following.
\begin{coro}{}\label{corostir}
Let $v=V/\log T$ be a real number, $V\in\RR^*$ bounded by an absolute constant. Let $h,k\leq T^\omega$, $0<\omega<1/2$ be positive integers. Then
\begin{align*}
&\frac{1}{\sqrt{hk}}\int_{T}^{2T}L(1/2+it,\chi)L(1/2-it-iv,\overline \chi)\left(\frac{h}k\right)^{-it}\left(\frac{2\pi}t\right)^{iv/2}\diff t=2T\frac{\chi\left(\frac{k}{(h,k)}\right)\overline\chi\left(\frac{h}{(h,k)}\right)}{hk}(h,k)\times 
\\&\left(L(1+iv,\chi_0)\left(\frac{q^2T(h,k)^2}{2\pi k^2}\right)^{iv/2}\frac{2^{1+iv/2}-1}{2+iv}+L(1-iv,\chi_0)\left(\frac{T(h,k)^2}{2\pi h^2}\right)^{-iv/2}\frac{2^{1-iv/2}-1}{2-iv}\right)
\\&+\tilde R(T;h,k;v).
\end{align*} 
Here, $\tilde R(T;h,k;v)$ is an error term satisfying the condition \eqref{condition}.
\end{coro}\vspace{-0.5cm}
\begin{proo}{}
First, note that $\left(|g(T)|+\int_1^T|g'(t)|\diff t\right)\ll \log T\ll T^{o(1)}$. Thus \nameref{coroipp} gives
\begin{align*}
&\frac{1}{\sqrt{hk}}\int_{1}^{T}L(1/2+it+iv,\chi)L(1/2-it,\overline \chi)\left(\frac{k}h\right)^{-it}\left(\frac{2\pi}t\right)^{-iv/2}\diff t=\frac{\chi\left(h/(h,k)\right)\overline\chi\left(k/(h,k)\right)}{hk}\times
\\&(h,k)\int_1^{T}\left(L(1-iv,\chi_0)\left(\frac{qt(h,k)}{2\pi k}\right)^{-iv}+\frac{L(1+iv,\chi_0)(h,k)^{iv}}{h^{iv}}\right)\left(\frac{2\pi}t\right)^{-iv/2}\diff t+R(T;k,h;iv,0;g)
\\&=2\frac{\chi\left(h/(h,k)\right)\overline\chi\left(k/(h,k)\right)}{hk}(h,k)\left(L(1-iv,\chi_0)\Biggl(\frac{(q(h,k))^2}{2\pi k^2}\right)^{-iv/2}\frac{T^{-iv/2+1}-1}{2-iv}
\\&\qquad +L(1+iv,\chi_0)\left(\frac{(h,k)^2}{2\pi h^2}\right)^{iv/2}\frac{T^{iv/2+1}-1}{2+iv}\Biggr)+R(T;k,h;iv,0;g).
\end{align*} 
Taking the conjugate of this, one finds the expected result.
\end{proo}

\section{Proof of \eqref{estiM}}
We start by the easy part of the proof of \nameref{prop:esti_inte}: the estimate on $M_j(t,H)$. 
\begin{proo}{of \eqref{estiM}} 
We fix $1\leq j\leq N$, we write $\chi$ for $\chi_j$, and we denote by $L$ the Dirichlet $L$-function associated to $\chi$. We also let $\eta=\eta_j$ be the corresponding mollifying function, and we write $M(t,H)$ instead of $M_j(t,H)$. 

We first apply Cauchy's residue theorem to the rectangle $1/2+i(t+H)$, $3/2+i(t+H)$, $3/2+it$, $1/2+it$: since $L(s)\eta(s)-1$ has no pole in that rectangle, we have
\begin{align*}
M(t,H)&=\int_t^{t+H}(L(3/2+iu)\eta^2(3/2+iu)-1)\diff u+\int_{1/2}^{3/2}(L(\sigma+it)\eta^2(\sigma+it)-1)\diff \sigma
\\&-\int_{1/2}^{3/2}(L(\sigma+i(t+H))\eta^2(\sigma+i(t+H))-1)\diff \sigma.
\end{align*}
Therefore, by the Cauchy-Schwarz inequality, we obtain
\begin{align*}
\numberthis\label{borneM}|M(t,H)|^2&\leq3\Biggl(\left|\int_t^{t+H}(L(3/2+iu)\eta^2(3/2+iu)-1)\diff u\right|^2+\left|\int_{1/2}^{3/2}(L(\sigma+it)\eta^2(\sigma+it)-1)\diff \sigma\right|^2
\\&+\left|\int_{1/2}^{3/2}(L(\sigma+i(t+H))\eta^2(\sigma+i(t+H))-1)\diff \sigma\right|^2\Biggr).
\end{align*}
We start by the first integral. By property of the Dirichlet convolution product, we know that 
\[\sum_{abc=n}\chi(a)\alpha(b)\alpha(c)=\delta_1(n),\]
where $\delta_1(n)=1$ if $n=1$, and $\delta_1(n)=0$ otherwise. Therefore, if $t\leq u\leq t+h$, the integrand of the first integral of \eqref{borneM} is 
\[\sum_{n,m,r=1}^\infty\frac{\chi(n)\beta(m)\beta(r)}{(nmr)^{3/2+iu}}-1=\sum_{n=1}^\infty\frac{\sum_{abc=n}\chi(a)\beta(b)\beta(c)}{n^{3/2+iu}}-1=\sum_{n\geq \xi^\theta}\frac{\sum_{abc=n}\chi(a)\beta(b)\beta(c)}{n^{3/2+iu}}.\]
The final equality relies on the fact that $\beta(k)=\alpha(k)$ if $k<\xi^\theta$. Moreover, 
\[\left|\sum_{abc=n}\chi(a)\beta(b)\beta(c)\right|\leq\sum_{abc=n}|\alpha(b)\alpha(c)|\leq \sum_{abc=n}1=\tau_3(n).\]
Note that $\tau_3(n)\ll n^{o(1)}$. Consequently, the first squared integral of \eqref{borneM} is 
\[\numberthis\label{first-esti}\ll \left(\sum_{n\geq  \xi^\theta}\frac{\tau_3(n)}{n^{3/2}\log n}\right)^2\ll \xi^{-\theta/2+o(1)}.\]
Next, a well known mean value estimate of Montgomery and Vaughan (see Corollary 3 of \cite{montgomery}) yields the following:
\[\numberthis\label{gab2}\int_T^{2T}\left|L(3/2+it)\eta^2(3/2+it)-1\right|^2\diff t\leq T\sum_{n\geq\xi^\theta}\frac{\tau_3(n)^2}{n^{3}}+\O(1)\lesssim \frac{T\xi^{-2\theta+o(1)}}2.\]
Now we turn to the second integral. From \nameref{lemme-mollifer}, the Cauchy-Schwarz inequality reveals that
\[\numberthis\label{gab1}\int_T^{2T}\left|L(1/2+it)\eta^2(1/2+it)-1\right|^2\diff t\leq (\frak C_3+1+2\sqrt{\frak C_3} )T.\]
Thanks to Gabriel's convexity theorem (see Theorem 2 of \cite{gabriel}), we use \eqref{gab1} and \eqref{gab2} to deduce that for all $1/2\leq \sigma\leq 3/2$:
\begin{align*}
\int_T^{2T}|L(\sigma+it)\eta^2(\sigma+it)-1|^2\diff t&\lesssim ((\frak C_3+1+2\sqrt{\frak C_3} )T)^{3/2-\sigma}\left(\frac{T\xi^{-2\theta+o(1)}}{2}\right)^{\sigma-1/2}
\\&\lesssim T \frac{(\frak C_3+1+2\sqrt{\frak C_3} )}{(2(\frak C_3+1+2\sqrt{\frak C_3} ))^{\sigma-1/2}}\xi^{-(\sigma-1/2)(2\theta+o(1))}.
\end{align*}
Therefore, by first applying the Cauchy-Schwarz inequality and then switching the integrals we find 
\begin{align*}
\int_T^{2T}&\left|\int_{1/2}^{3/2}(L(\sigma+it)\eta^2(\sigma+it)-1)\diff \sigma\right|^2\diff t
\\&\leq \int_T^{2T}\left[\int_{1/2}^{3/2} \xi^{-(\sigma_1-1/2)\theta}\diff \sigma_1\int_{1/2}^{3/2} \xi^{(\sigma_2-1/2)\theta}|L(\sigma+it)\eta^2(\sigma_2+it)-1|^2\diff\sigma_2\right]\diff t
\\&\lesssim T (\frak C_3+1+2\sqrt{\frak C_3} )\left[-\frac{\xi^{-(\sigma_1-1/2)\theta}}{\theta\log\xi}\right]_{\sigma_1=1/2}^{3/2}\int_{1/2}^{3/2}\frac{\xi^{-(\sigma_2-1/2)(\theta+o(1))}}{(2(\frak C_3+1+2\sqrt{\frak C_3} ))^{\sigma_2-1/2}}\diff \sigma_2
\\&\lesssim \frac{(\frak C_3+1+2\sqrt{\frak C_3} )T}{(\theta\log \xi)^2}\lesssim\frac{(\frak C_3+1+2\sqrt{\frak C_3} )}{(\theta\varkappa)^2} \frac{T}{(\log T)^2},
\end{align*}
since $\int_{1/2}^{3/2}(\xi^{\theta+o(1)}/C)^{-(\sigma-1/2)}\diff \sigma=(1+o(1))/(\theta\log \xi)$ for any constant $C\neq 0$. The same estimate holds for the third integral of \eqref{borneM}. Therefore, combining these estimates with \eqref{first-esti}, we get that 
\[\int_T^{2T}|M(t,H)|^2\diff t\lesssim 3\left(T\xi^{-\theta/2+o(1)}+\frac{2(\frak C_3+1+2\sqrt{\frak C_3} )}{(\theta\varkappa)^2} \frac{T}{(\log T)^2}\right)\lesssim\frac{6(\frak C_3+1+2\sqrt{\frak C_3} )}{(\theta\varkappa)^2} \frac{T}{(\log T)^2}.\]
\end{proo}

\section{Proof of \nameref{lemme-mollifer}}
We may now prove \nameref{lemme-mollifer}. Again, we fix $1\leq j\leq N$, we write $\chi$ for $\chi_j$, and we denote by $L$ the Dirichlet $L$-function associated to $\chi$. We also let $\eta=\eta_j$ be the corresponding mollifying function.
\begin{proo}{of \nameref{lemme-mollifer}}
We suppose that $\varkappa<2/9$ and we fix $\varepsilon>0$ small enough. We let $(\gamma(n))$ be the coefficients such that 
\[\eta(s)^2=\sum_{n}\frac{\gamma(n)\chi(n)}{n^s}.\]
Note that $\gamma(n)=\mathcal M\tau_{-1/2}\star\mathcal M\tau_{-1/2}(n)$, where $\star$ stands for the Dirichlet convolution, and hence is real and is supported on integers $n\leq \xi^2$. By definition, we have that $|\gamma(n)|\leq \tau_{1/2}\star\tau_{1/2}(n)=1$. Also note that, if $\chi_0$ is the principal character $\mod q$, then for any integers $h,k$:
\[\numberthis\label{charprin}\chi(h)\overline{\chi}(k)\overline{\chi}\left(\frac{h}{(h,k)}\right)\chi\left(\frac{k}{(h,k)}\right)=\chi_0\left(\frac{hk}{(h,k)}\right)=\chi_0(hk).\]

Then, by \nameref{coroegal} with $\omega=2\varkappa$, we get
\begin{align*}
&\int_T^{2T}|L(1/2+it)\eta^2(1/2+it)|^2\diff t=\sum_{h,k}\frac{\gamma(h)\gamma(k)\chi(h)\overline\chi(k)}{\sqrt{hk}}\int_T^{2T}|L(1/2+it)|^2\left(\frac hk\right)^{-it}\diff t
\\&=\sum_{h,k}\frac{\chi_0(hk)\gamma(h)\gamma(k)(h,k)}{hk}\frac{\phi(q)}{q}T\left(\log\left(\frac{T(h,k)^2}{hk}\right)+\frak Q(q,0)+2\log 2\right)
\\&\qquad+\O((T^{1/2+2\varkappa}+T^{1/3+3\varkappa})T^\varepsilon)
\\\numberthis\label{eqlemme4}&=\frac{\phi(q)}q\Biggl[\Sigma_1 T(\log T+\frak Q(q,0)+2\log 2)-\Sigma_2 T\Biggr]+o(T),
\end{align*}
where 
\[\Sigma_1:=\sum_{h,k}\frac{\chi_0(hk)\gamma(h)\gamma(k)(h,k)}{hk},\]
\[\Sigma_2:=\sum_{h,k}\frac{\chi_0(hk)\gamma(h)\gamma(k)(h,k)}{hk}\log\left(\frac{hk}{(h,k)^2}\right).\]
Now we only have to bound these two sums, as in \cite[24.51]{iwaniec}. We have
\begin{align*}
\Sigma_1&=\sum_{d}d^{-1}\sum_{\substack{(h,k)=1}}\frac{\chi_0(d^2hk)\gamma(dh)\gamma(dk)}{hk}=\sum_{d}d^{-1}\sum_{\delta}\mu(\delta)\sum_{\substack{h,k}}\frac{\chi_0(d^2\delta ^2hk)\gamma(d\delta h)\gamma(d\delta k)}{\delta^2 hk}
\\&\numberthis\label{eqAdfin}=\sum_{m}m^{-1}\sum_{\delta|m}\frac{\mu(\delta)}{\delta}\left(\sum_{\substack{h}}\frac{\chi_0(mh)\gamma(mh)}{h}\right)^2=\sum_{d}\phi(d)A_d^2,
\end{align*}
where 
\[A_d=\frac{1}{d}\sum_{h}\frac{\chi_0(dh)\gamma(dh)}h=\sum_{\substack{h\equiv 0\mod d}}\frac{\chi_0(h)\gamma(h)}h.\]
Because of the presence of $\chi_0$, we may restrict our attention to $d$ co-prime to $q$. For any positive integers $h_1,h_2$ such that $h_1h_2=h$, we let $h_1=\delta_1 n$, $h_2=\delta_2 m$, $\delta_1\delta_2|d^\infty$, and $(nm,d)=1$. Note that since $\delta_1\delta_2|d^\infty$ and $(q,d)=1$, we have $(\delta_1\delta_2,q)=1$. With this decomposition we may write
\begin{align*}
A_d&=\sum_{\substack{h\equiv 0\mod d}}\sum_{h_1h_2=h}\frac{\mathcal M(h_1)\chi_0(h_1)\tau_{-1/2}(h_1)\mathcal M(h_2)\chi_0(h_2)\tau_{-1/2}(h_2)}{h_1h_2}
\\&=\sum_{\substack{\delta_1\delta_2|d^\infty\\ d|\delta_1\delta_2}}\sum_{\substack{n,m\geq 1\\(nm,d)=1}}\frac{\chi_0(n)\mathcal M(\delta_1n)\tau_{-1/2}(\delta_1n)\chi_0(m)\mathcal M(\delta_2m)\tau_{-1/2}(\delta_2m)}{\delta_1\delta_2nm}
\\&\numberthis\label{eqAddelta}=\sum_{\substack{\delta_1\delta_2|d^\infty\\ d|\delta_1\delta_2}}\frac{B_d(\delta_1)B_d(\delta_2)}{\delta_1\delta_2}
\end{align*}
where we have defined $B_d(\delta):=\sum_{(n,qd)=1}\frac{\mathcal M(n\delta)\tau_{-1/2}(n\delta)}{n}$, for simplicity. We may use \nameref{lemmeA} to get that
\begin{align*}
|B_d(\delta)|\numberthis\label{defC'}\lesssim\frak C_5\tau_{1/2}(\delta)\sqrt{\frac{dq}{(\log T)\phi(d)\phi(q)}},
\end{align*}
where 
\[\frak C_5:=\frac1{1-\theta}\frac{e^{\varrho_\theta}+e^{\varrho_\theta \theta}}{2 \sqrt{\pi\varkappa\varrho_\theta}}\frac{\Gamma(1/4)}{\Gamma(3/4)}\]
and where $\varrho_\theta$ is the only positive solution of $-1+2\theta x+e^{x(1-\theta)}(-1+2x)=0$.

Since 
\[\sum_{r|d^\infty, \, d|r}\frac{\tau_{1/2}\star\tau_{1/2}(r)}{r}=\sum_{r|d^\infty}\frac{1}{dr}=\frac1d\prod_{p|d}(1-p^{-1})^{-1}=\frac1{\phi(d)},\]
putting \eqref{defC'} in \eqref{eqAddelta}, we find
\begin{align*}
|A_d|\lesssim \sum_{\substack{\delta_1\delta_2|d^\infty\\ d|\delta_1\delta_2}}\frak{C}_5^2\frac{dq}{(\log T)\phi(d)\phi(q)}\frac{\tau_{1/2}(\delta_1)\tau_{1/2}(\delta_2)}{\delta_1\delta_2}&\leq \frak{C}_5^2\frac{dq}{(\log T)\phi(d)\phi(q)}\sum_{r|d^\infty, \, d|r}\frac{\tau_{1/2}\star\tau_{1/2}(r)}{r}
\\\numberthis\label{eqadfinal}&=\frak{C}_5^2\frac{dq}{(\log T)\phi(d)^2\phi(q)},
\end{align*}
Finally, since $A_d$ is zero if $(d,q)\neq 1$, inserting \eqref{eqadfinal} in \eqref{eqAdfin} yields
\[\numberthis\label{sig1}|\Sigma_1|\lesssim\frak{C}_5^4\frac{q^2}{(\log T)^2\phi(q)^2}\sum_{d\leq T^{2\varkappa}}\frac{d^2}{\phi(d)^3}\chi_0(d).\]
Observe that for $p$ a prime and $m\geq 1$ an integer, we have
\[\frac{p^{2m}}{\phi(p^m)^3}=\frac{1}{p^m}\left(\frac{p}{p-1}\right)^3=\frac{1}{p^m}+\frac{1}{p^{m-1}}\left(\frac1p\left(\frac{p}{p-1}\right)^3-\frac1p\right)=\frac{1}{p^m}+\frac{1}{p^{m-1}}\left(\frac{3p^2-3p+1}{p^4-3p^3+3p^2-p}\right).\]
Thus using multiplicativity we find that for $d\geq 1$ an integer, we have
\[\frac{d^2}{\phi(d)^3}=inv\star\mu^2\psi(d),\]
where $\psi$ is the totally multiplicative function defined on primes by $p\mapsto (3p^2-3p+1)/(p^4-3p^3+3p^2-p)$ and $inv(d)=1/d$. Thus
\begin{align*}
\sum_{d\leq T^{2\varkappa}}\frac{d^2}{\phi(d)^3}\chi_0(d)&=\sum_{d\leq T^{2\varkappa}}\mu^2(d)\psi(d)\chi_0(d)\sum_{n\leq  T^{2\varkappa}/d}\frac{\chi_0(n)}n\leq \sum_{n\leq T^{2\varkappa}}\frac{\chi_0(n)}n\sum_{d=1}^{\infty}\mu^2(d)\psi(d)
\\&\numberthis\label{eqprodeul}\lesssim 2\varkappa (\log T)\frac{\phi(q)}q\prod_{p}\left(1+\frac{3p^2-3p+1}{p^4-3p^3+3p^2-p}\right).
\end{align*}

Therefore, using \eqref{sig1}, we get
\[\numberthis\label{sigma1}|\Sigma_1|\lesssim 2\varkappa \frac{q}{\phi(q)}\frak{C}_5^4\prod_{p}\left(1+\frac{3p^2-3p+1}{p^4-3p^3+3p^2-p}\right) (\log T)^{-1}.\]

Similarly, we will show that $\Sigma_2\ll 1$ (see \cite[24.50]{iwaniec}). First, we have 
\[\Sigma_2=\sum_{d}\frac1d\sum_{(h,k)=1}\frac{\chi_0(d^2hk)\gamma(dh)\gamma(dk)}{hk}\log(hk).\]
We write $\log(hk)=\log(hkd^2)-2\log d$, and split $\Sigma_2$ in two according to this decomposition, say $\Sigma_2=\Sigma_2'-2\Sigma_2''$. Similarly to the case of $\Sigma_1$, we have 
\[\Sigma_2'=2\sum_d\phi(d)A_d\tilde A_d,\]
where
\[\tilde A_d=\sum_{d|n}\frac{\gamma(n)\chi_0(n)\log n}n=\sum_{d|n}\sum_{k|n}\Lambda(k)\frac{\gamma(n)\chi_0(n)}{n}=\sum_{k|n}\Lambda(k)\sum_{[d,k]|n}\frac{\gamma(n)\chi_0(n)}{n}=\sum_{k}\Lambda(k)A_{[d,k]}.\]
Here, $\Lambda$ is von Mangoldt's function, defined by $1\star\Lambda= \log$. Using \eqref{eqadfinal} and the fact that for integers $a,b\geq 1$ we have $\phi([a,b])\phi((a,b))=\phi(a)\phi(b)$, one gets 
\[\numberthis\label{bd}\tilde A_d\lesssim \frac{\frak{C}_5^2}{\log T} \frac{q}{\phi(q)}\sum_{k\leq T^{2\varkappa}}\frac{[d,k]\Lambda(k)}{\phi^2([d,k])}=\frac{\frak{C}_5^2dq}{\phi^2(d)\phi(q)\log T}\sum_{k\leq T^{2\varkappa}}\frac{k\Lambda(k)\phi^2((d,k))}{(d,k)\phi^2(k)}.\]
If $\nu_p$ is the $p$-adic valuation, $P>0$ is real number, and $d\leq T^{2\varkappa}$, then this last sum is 
\begin{align*}
&\leq\sum_{p|d}\sum_{p^m\leq T^{2\varkappa},m\geq 1} \frac{(\log p)p^{\min(\nu_p(d),m)}}{p^m}+\sum_{p\nmid d}\sum_{p^m\leq T^{2\varkappa},m\geq 1}\frac{\log p}{p^m}\left(\frac{p}{p-1}\right)^2
\\&\leq\sum_{p|d}\sum_{1\leq m\leq\nu_p(d)} \log p+\sum_{p|d}\sum_{p^m\leq T^{2\varkappa}, m\geq 1}\frac{\log p}{p^m}+\sum_{p\nmid d}\sum_{p^m\leq T^{2\varkappa},m\geq 1}\left(\frac{p}{p-1}\right)^2\frac{\log p}{p^m}
\\&\leq\sum_{p^e||d}\log(p^e)+\sum_{p}\sum_{p^m\leq T^{2\varkappa},m\geq 1}\left(\frac{p}{p-1}\right)^2\frac{\log p}{p^m}
\\&\leq \log d+4\sum_{p\leq P}\sum_{p^m\leq T^{2\varkappa},m\geq 1}\frac{\log p}{p^m}+\left(\frac{P}{P-1}\right)^2\sum_{P<p\leq T^{2\varkappa}}\sum_{p^m\leq T^{2\varkappa},m\geq 1}\frac{\log p}{p^m}
\\&\leq 2\varkappa\log T+4\log P\sum_{p\leq P}\frac{1}{p-1}+\left(1+\O\left(\frac 1P\right)\right)\sum_{k\leq T^{2\varkappa}}\frac{\Lambda(k)}{k}.
\end{align*}
Note that the last sum is $\log (T^{2\varkappa})+\O(1)$. Taking $P=\log T$ for example, one may inject this in \eqref{bd} to find that:
\[\tilde A_d\lesssim 4\varkappa\frac{\frak{C}_5^2dq}{\phi^2(d)\phi(q)}.\]

Thus, using \eqref{eqadfinal} and \eqref{eqprodeul}, we find
\[\numberthis\label{sigma'}\Sigma_2'\lesssim \frac{4\varkappa\frak{C}_5^4}{\log T}\frac{q^2}{\phi(q)^2}\sum_{d\leq T^{2\varkappa}}\frac{d^2}{\phi(d)^3}\chi_0(d)\lesssim 8\varkappa^2\frak{C}_5^4\frac{q}{\phi(q)}\prod_{p}\left(1+\frac{3p^2-3p+1}{p^4-3p^3+3p^2-p}\right).\]
As in \eqref{eqAdfin}, by successively making the variable changes $r=d/\delta$ and $m=r/k$, we find that
\begin{align*}
\Sigma_2''&=\sum_{d}A_d^2\sum_{\delta|d}\mu(\delta)\frac d\delta\log(d/\delta)=\sum_{d}A_d^2\sum_{r|d}\sum_{k|r}\Lambda(k)r\mu(d/r)
\\&=\sum_{d}A_d^2\sum_{k|d}\Lambda(k)\sum_{m|\frac dk}km\mu\left(\frac d{km}\right)=\sum_{d}A_d^2\sum_{k|d}\Lambda(k)k\phi(d/k).
\end{align*}
Using the fact that for integers $a,b\geq 1$, $\phi(a)\phi(b)\leq \phi(ab)$, we find that the inner sum above is
\begin{align*}
\leq\phi(d)\sum_{k|d}\Lambda(k)\frac{k}{\phi(k)}= \phi(d)\sum_{p^m|d,m\geq 1}\log p \frac{p}{p-1}\leq 2\phi(d)\sum_{p^m|d,m\geq 1}\log p &= 2\phi(d)\log d
\\&\leq 4\varkappa \phi(d)\log T.
\end{align*}
Therefore \eqref{eqadfinal} and \eqref{eqprodeul} yield
\[\Sigma_2''\lesssim 8\varkappa^2\frak{C}_5^4\frac{q}{\phi(q)} \prod_{p}\left(1+\frac{3p^2-3p+1}{p^4-3p^3+3p^2-p}\right).\]
Together with \eqref{sigma'}, \eqref{sigma1} and \eqref{eqlemme4}, this concludes the proof.
\end{proo}

\section{Proof of \eqref{estiI}}
Finally we turn to the last estimate needed to conclude the proof of \nameref{prop:esti_inte}. We fix $1\leq j\leq N$, we write $\chi$ for $\chi_j$, and we denote by $L$ the Dirichlet $L$-function associated to $\chi$. We also let $\eta=\eta_j$ be the corresponding mollifying function, and we write $I(t,H)$ instead of $I_j(t,H)$. 
\begin{proo}{of \eqref{estiI}}
Write $Y(u)=X(1/2+iu)|\eta^2(1/2+iu)|$ for convenience. Using the subconvexity bound $L(1/2+it,\chi)\ll t^{1/6+o(1)}$ for any $t>2$ and $\varepsilon>0$ (see Theorem 2 of \cite{kol}) and since $|\gamma(n)|\leq 1$, one gets
\begin{align*}
\int_T^{2T}&\left|\int_t^{t+H}Y(u)\diff u\right|^2\diff t=\int_0^H\int_0^H\int_T^{2T}Y(u+t)\overline{Y(v+t)}\diff t\diff u\diff v
\\&=\int_0^H\int_0^H\int_{T+u}^{2T+u}Y(t)\overline{Y(t+v-u)}\diff t\diff u\diff v
\\&=\int_0^H\int_0^H\int_{T}^{2T}Y(t)\overline{Y(t-u+v)}\diff t\diff u\diff v+\O\left(T^{1/3+o(1)}\left(\sum_{n\leq \xi^2}\frac{|\gamma(n)|}{\sqrt{n}}\right)^2\right)
\\&=\int_0^H\int_{-u}^{H-u}\int_{T}^{2T}Y(t)\overline{Y(t+w)}\diff t\diff w\diff u+\O(T^{2\varkappa+1/3+o(1)})
\\&=\numberthis\label{eqintu}\int_{-H}^H(H-|v|)\int_T^{2T}Y(t)\overline{Y(t+v)}\diff t\diff v+\O(T^{2\varkappa+1/3+o(1)}).
\end{align*}
The last equality is obtained by switching the first two integrals. By Stirling's formula for the $\Gamma$-function (see Corollary II.0.13 of \cite{tenenbaum}), we may write that for $|v|\leq H$:
\[e^{i\vartheta(1/2+it)}e^{-i\vartheta(1/2+it+iv)}=\left(\frac{2\pi}{t}\right)^{iv/2}\left(1+\O\left(\frac{1}t\right)\right).\]
We let $\tilde{\gamma}(n):=\mathcal{M}(n)\tau_{-1/2}(n)$ for simplicity. For $t\in(T,2T)$ and $|v|\leq H$, we then deduce -again by the subconvexity bound- that
\begin{align*}
&Y(t)\overline{Y(t+v)}=q^{-iv/2}L(1/2+it,\chi)L(1/2-it-iv,\overline\chi)|\eta^2(1/2+it)\eta^2(1/2+it+iv)|\left(\frac{2\pi}{t}\right)^{iv/2}
\\&\qquad +\O(T^{2\varkappa-2/3+o(1)})
\\&=\sum_{h,k,h',k'}L(1/2+it,\chi)L(1/2-it-iv,\overline\chi)\frac{\beta(h)\overline\beta(k)\beta(h')\overline\beta(k')}{\sqrt{hkh'k'}}\left(\frac{kk'}{hh'}\right)^{it}\left(\frac{2\pi k'^2}{tqh'^2 }\right)^{iv/2}
\\&\qquad +\O(T^{2\varkappa-2/3+o(1)})
\\&=\sum_{h,k,h',k'}L(1/2+it,\chi)L(1/2-it-iv,\overline\chi)\frac{\tilde\gamma(h)\tilde\gamma(k)\tilde\gamma(h')\tilde\gamma(k')}{\sqrt{hkh'k'}}\chi(hh')\overline\chi(kk')\left(\frac{kk'}{hh'}\right)^{it}\left(\frac{2\pi k'^2}{tqh'^2 }\right)^{iv/2}
\\&\qquad +\O(T^{2\varkappa-2/3+o(1)}).
\end{align*}
Note that over our range of summation, we have $hh',kk'\leq T^{2\varkappa}$, and recall \eqref{charprin}. Integrating this over $t$, our \nameref{corostir} applied with $\omega=2\varkappa$ yields 
\begin{align*}
&\int_T^{2T}Y(t)\overline{Y(t+v)}\diff t=2T\sum_{h,k,h',k'}\left(\frac{k'^2}{qh'^2 }\right)^{iv/2}\frac{\chi_0(kk'hh')\tilde\gamma(h)\tilde\gamma(k)\tilde\gamma(h')\tilde\gamma(k')}{hkh'k'}(hh',kk')\times
\\&\Biggl(L(1+iv,\chi_0)\left(\frac{q^2T(hh',kk')^2}{2\pi (kk')^2}\right)^{iv/2}\frac{2^{1+iv/2}-1}{2+iv}+L(1-iv,\chi_0)\left(\frac{T(hh',kk')^2}{2\pi (hh')^2}\right)^{-iv/2}\frac{2^{1-iv/2}-1}{2-iv}\Biggr)
\\&\numberthis\label{termerr1}\qquad+\O\left(\sum_{h,h',k,k'}\left(\frac{k'^2}{h'^2 }\right)^{iv/2}\chi(hh')\overline\chi(kk')\tilde R(T;hh',kk';iv)\right)+\O(T^{2\varkappa+1/3+o(1)}).
\end{align*}
By property of $\tilde R$ and by letting $K'=kk'$ and $H'=hh'$, the first error term is 
\begin{align*}
&\ll \sum_{h,k\leq\xi}\left(\frac hk\right)^{iv/2}\sum_{K',H'\leq \xi^2}\left(\frac{K'^2}{H'^2 }\right)^{iv/2}\chi(H')\overline\chi(K')\frak 1_{k|K'}\frak 1_{h|H'}\tilde R(T;H',K';iv)
\\&\ll(T^{1/2+2\varkappa}+T^{1/3+3\varkappa})T^{o(1)}\sum_{h,k\leq \xi}1 \ll (T^{1/2+4\varkappa}+T^{1/3+5\varkappa})T^{o(1)},
\end{align*}
where $\frak 1_{n|m}=1$ if $n|m$ and $\frak 1_{n|m}=0$ otherwise.

Because $\varkappa<1/8$, the error term of \eqref{termerr1} is $\ll T^{1-\epsilon}$, for any small $\epsilon>0$. Using the symmetric roles of the variables, one gets
\begin{align*}
&\int_T^{2T}Y(t)\overline{Y(t+v)}\diff t=4T\Re\Biggl[\sum_{h,k,h',k'}\left(\frac{k'}{h'}\right)^{iv}\frac{\chi_0(hh'kk')\tilde\gamma(h)\tilde\gamma(k)\tilde\gamma(h')\tilde\gamma(k')}{hkh'k'}(hh',kk')
\\&\times L(1+iv,\chi_0)\left(\frac{qT(hh',kk')^2}{2\pi (kk')^2}\right)^{iv/2}\frac{2^{1+iv/2}-1}{2+iv}\Biggr]+\O(T^{1-\epsilon})
\\&\qquad\qquad\qquad\qquad \quad=4T\Re\left[K(v)L(1+iv,\chi_0)\frac{2^{1+iv}-1}{2+iv}\right]+\O(T^{1-\epsilon}),
\end{align*}
where
\[K(v):=\sum_{h,k,h',k'}\frac{\chi_0(hh'kk')\tilde\gamma(h)\tilde\gamma(k)\tilde\gamma(h')\tilde\gamma(k')}{hkh'k'}(hh',kk')\left(\frac{(hh',kk')}{h'k}\sqrt{\frac{qT}{2\pi}}\right)^{iv}.\]
Since $L(s,\chi_0)$ only has a pole of order $1$ at $s=1$, of residue $\phi(q)/q$, we may expend $L$ as a Laurent series to find that
\[\frac{L(1+iv,\chi_0)}{2+iv}(2^{1+iv/2}-1)=\frac{\phi(q)}{2iqv}+\O(1),\]
and hence we may use \eqref{eqintu} to deduce that
\begin{align*}
\int_T^{2T} |I(t,H)|^2\diff t&=2T\frac{\phi(q)}q\int_{-H}^H(H-|v|)\Re\left[\frac{K(v)}{iv}\right]\diff v+\O\left(T\int_{-H}^H(H-|v|)\left|K(v)\right|\diff v\right)
\\\numberthis\label{expressionI}&\qquad+\O(T^{1-\epsilon}).
\end{align*}
We fix $|v|\leq H$. We will estimate $K(v)$ in a similar manner to the one used to estimate $\Sigma_1$, in the proof of \nameref{lemme-mollifer}. We will also add some arguments of arithmetical nature that may be found in \cite{selberg_zeta}. We will divide the rest of the proof into three parts. First, we will show that $K(v)=\O(1/\log T)$, thus showing that the error term is $\O(T/(\log T)^3)$. Then, we will prove that $\Im(K(v)/v)$ is $\O(1)$ if $|v|\leq 1/\log \xi$. Finally we will bound it for the remaining $v$'s. 

We let, for $z$ a complex number, $\phi_z(n):=\sum_{m|n}\mu(m)\left(\frac{n}m\right)^{1+z}$. By Möbius inversion, it is easy to see that if $n\geq 1$ is an integer, then 
\[\numberthis\label{eqmobius}n^{1+z}=\sum_{m|n}\phi_z(m).\]
Note that for any real number $|v|\leq H$, 
\[\numberthis\label{majIm2}|\phi_{iv}(d)|\leq d\sum_{\rho|d}\frac1\rho.\]
Again we put $\mathcal T=qT/(2\pi)$. Using \eqref{eqmobius}, we write the following:
\begin{align*}
K(v)&=\mathcal T^{iv/2}\sum_{h,h',k,k'}\frac{\chi_0(hh'kk')\tilde{\gamma}(h)\tilde{\gamma}(h')\tilde{\gamma}(k)\tilde{\gamma}(k')}{hh'kk'}\frac{(hh',kk')^{1+iv}}{(h'k)^{iv}}
\\&=\mathcal T^{iv/2}\sum_{d\leq T^{2\varkappa}}\phi_{iv}(d)\sum_{\substack{h,k,h',k'\\d|hh'\\d|kk'}}\frac{\chi_0(hh'kk')\tilde{\gamma}(h)\tilde{\gamma}(h')\tilde{\gamma}(k)\tilde{\gamma}(k')}{h'^{1+iv}k^{1+iv}hk'}
\\\numberthis\label{eqExpK}&=\mathcal T^{iv/2}\sum_{d\leq T^{2\varkappa}}\phi_{iv}(d)\left(\sum_{\substack{h,h'\\d|hh'}}\frac{\chi_0(hh')\tilde{\gamma}(h)\tilde{\gamma}(h')}{h'h^{1+iv}}\right)^2.
\end{align*}
We denote the inner sum by $A'_d(v)$, similarly to the quantity $A_d$ introduced in the proof of \nameref{lemme-mollifer}. Again, we may restrict our attention to $d$ that are co-prime to $q$. We let $\delta,\delta'|d^\infty$ (and hence $(\delta\delta',q)=1$) be such that $h=\delta k$, $h'=\delta' k'$ for some $k,k'$ such that $(kk',d)=1$. With these notations, we have
\[\numberthis\label{eqa'}A'_d(v)=\sum_{\substack{d|\delta\delta'\\\delta\delta'|d^\infty}}\frac1{\delta^{1+iv}\delta'}\left(\sum_{(k,qd)=1}\frac{\mathcal{M}(k\delta)\tau_{-1/2}(k\delta)}{k^{1+iv}}\right)\left(\sum_{(k,qd)=1}\frac{\mathcal{M}(k\delta')\tau_{-1/2}(k\delta')}{k}\right).\]
The sum inside the first pair of parenthesis is $S_d(\delta,v)$, and the one in the second pair of parenthesis is simply equal to $B_d(\delta')$, these quantities being introduced in \nameref{lemmeA} and \nameref{lemme-mollifer} respectively. Then \eqref{lemmeA} leads to 
\begin{align*}
|S_d(\delta,v)|&\lesssim\frak{C}_6(v\log T)\tau_{1/2}(\delta)\sqrt{\frac{dq}{(\log T)\phi(d)\phi(q)}},
\end{align*}
where  
\[\frak{C}_6(v\log T):=\frac{1}{1-\theta}\frac{e^{\rho(v\log T,\theta)}+e^{\rho(v\log T,\theta) \theta}}{2\sqrt{\pi\varkappa\rho(v\log T,\theta)}}\left(\sqrt{\frac{v\log T}{\rho(v\log T,\theta)}}\sqrt{\pi\varkappa}+\frac{\Gamma(1/4)}{\Gamma(3/4)}\right).\]
By \eqref{defC'}, we also have
\[|B_d(\delta')|\lesssim\frak{C}_5\tau_{1/2}(\delta')\sqrt{\frac{dq}{(\log T)\phi(d)\phi(q)}}.\]
These inequalities lead to (see the similar computation \eqref{eqadfinal}):
\[\numberthis\label{eqbesoin}|A'_d(v)|\lesssim \frak{C}_5\frak C_6(v\log T)\sum_{\substack{d|\delta\delta'\\\delta\delta'|d^\infty}}\frac{\tau_{1/2}(\delta)\tau_{1/2}(\delta')}{\delta\delta'}\frac{dq}{(\log T)\phi(d)\phi(q)}\leq \frak{C}_5\frak C_6(v\log T)\frac{dq}{(\log T)\phi(d)^2\phi(q)}.\]
Therefore, by \eqref{majIm2} and \eqref{eqExpK}, we have 
\begin{align*}
K(v)&\ll \sum_{d\leq T^{2\varkappa}} d\sum_{\rho|d}\frac 1\rho \left(\frac{d}{(\log T)\phi(d)^2}\right)^2\ll\frac{1}{(\log T)^2}\sum_{d\leq T}\frac{d^3\sum_{\rho|d}\frac 1\rho}{\phi(d)^4}
\\&\ll\frac{1}{(\log T)^2}\sum_{d\leq T} \frac{1}d\prod_{p|d}(1-1/p)^{-5}\ll\frac{1}{(\log T)^2}\sum_{d\leq T} \frac{1}d\prod_{p|d}(1+1/\sqrt p)
\\&\ll\frac{1}{(\log T)^2} \sum_{d\leq T}d^{-1}\sum_{\rho|d}\rho^{-1/2}\leq\frac{1}{(\log T)^2} \zeta(3/2)\sum_{d\leq T}\frac1d\ll \frac{1}{\log T}.
\end{align*}
This proves that \eqref{expressionI} can be written as 
\begin{align*}
\int_T^{2T} |I(t,H)|^2\diff t&=2T\frac{\phi(q)}q\int_{-H}^H(H-|v|)\Im\left[\frac{K(v)}{v}\right]\diff v+\O\left(\frac{T}{(\log T)^3}\right)
\\\numberthis\label{eqIntfacile}&=4T\frac{\phi(q)}q\int_{0}^H(H-v)\Im\left[\frac{K(v)}{v}\right]\diff v+\O\left(\frac{T}{(\log T)^3}\right).
\end{align*}

Now we are interested about this integral over the interval $0<v<1/\log \xi$, and we fix such a $v$. Observe that for $a,b\in\CC$, we have 
\[|\Im(a^2b)|=|\Im(a^2)\Re(b)+\Re(a^2)\Im(b)|\leq |b||\Im(a^2)|+|a|^2|\Im(b)|\leq 2|a||b||\Im(a)|+|a|^2|\Im(b)|.\]
Then
\begin{align*}
\left|\Im\left(\frac{K(v)}{v}\right)\right|&\leq\sum_{d\leq T^{2\varkappa}}\left|\Im\left(\mathcal T^{iv/2}\frac{\phi_{iv}(d)}vA'_d(v)^2\right)\right|\chi_0(d)
\\\numberthis\label{ineImaPart}&\leq \sum_{d\leq T^{2\varkappa}}\left[\left|\Im\left(\frac{\mathcal T^{iv/2}\phi_{iv}(d)}{v}\right)\right||A'_d(v)|^2+2|\phi_{iv}(d)||A_d'(v)|\left|\Im\left(\frac{A_d'(v)}v\right)\right|\right]\chi_0(d).
\end{align*}

Since $\sin(x)\leq x$ for any $x>0$, and by definition of $\phi_{iv}$, we have for $d\leq T^{2\varkappa}$, $(d,q)=1$:
\begin{align*}
\numberthis\label{majIm}\left|\Im\left(\frac{\mathcal T^{iv/2}\phi_{iv}(d)}{v}\right)\right|&=d\left|\sum_{\rho|d}\frac{\mu(\rho)}{\rho}\frac{\sin(v\log(\sqrt{\mathcal T} d/\rho))}v\right|\lesssim d\left(\frac12+2\varkappa\right)\log T\sum_{\rho|d}\frac 1\rho.
\end{align*}
Also, 
\begin{align*}
&\Im\left(\frac{A_d'(v)}v\right)=\sum_{\substack{h,h'\\d|hh'}}\frac{\chi_0(hh')\tilde{\gamma}(h)\tilde{\gamma}(h')}{vh'h}\sin(-v\log h)
\\&=\sum_{\substack{d|\delta\delta'\\\delta\delta'|d^\infty}}\frac1{\delta\delta'}\left(\sum_{(k,qd)=1}\frac{\mathcal{M}(k\delta)\tau_{-1/2}(k\delta)}{kv}\sin(-v\log(k\delta))\right)\left(\sum_{(k,qd)=1}\frac{\mathcal{M}(k\delta')\tau_{-1/2}(k\delta')}{k}\right),
\end{align*}
and hence by \nameref{coro-molli2} and \eqref{defC'}
\begin{align*}
\left|\Im\left(\frac{A_d'(v)}v\right)\right|&\lesssim\sum_{\substack{d|\delta\delta'\\\delta\delta'|d^\infty}}\frac{\tau_{1/2}(\delta)}{\delta\delta'}\left(\frak C_4\sqrt{\frac{dq}{\phi(d)\phi(q)}}\sqrt{\log \xi}+\O\left(\prod_{p|d}(1+p^{-3/4})(\log T)^{1/3}\right)\right)
\\&\qquad\times \frak{C}_5\tau_{1/2}(\delta')\sqrt{\frac{dq}{(\log T)\phi(d)\phi(q)}}
\\&\numberthis\label{EqImAd}\leq \frak C_4\frak C_5\sqrt\varkappa \frac{dq}{\phi(d)^2\phi(q)}+\O\left(\frac{\sqrt d}{\phi(d)^{3/2}(\log T)^{1/6}}\prod_{p|d}(1+p^{-3/4})\right).
\end{align*}
The last inequality is justified as in \eqref{eqadfinal}. 

Collecting \eqref{majIm}, \eqref{eqbesoin}, \eqref{majIm2}, \eqref{EqImAd} and injecting them in \eqref{ineImaPart}, we get
\begin{align*}
\Im\left(\frac{K(v)}v\right)&\lesssim \frac{1}{\log T} \left[\left(\frac12+2\varkappa\right)(\frak C_5\frak C_6(v\log T))^2+2\frak C_4\frak C_5^2\frak C_6(v\log T)\sqrt\varkappa \right]\frac{q^2}{\phi(q)^2}\sum_{d\leq \xi^2}\frac{d^3\chi_0(d)}{\phi(d)^4}\sum_{\rho|d}\frac1\rho
\\&\numberthis\label{eq-Im_Kv/v}\qquad+\O\left(\frac1{(\log T)^{7/6}}\sum_{d\leq \xi^2}\frac{d^{5/2}\sum_{\rho|d}\frac1\rho}{\phi(d)^{7/2}}\prod_{p|d}(1+p^{-3/4})\right).
\end{align*}
To estimate the double sum in the main term, we let $f(d):=\left(\frac{d}{\phi(d)}\right)^4\sum_{\rho|d}\frac{1}{\rho}$, which is multiplicative, and we let $g=\mu\star f$, which is also multiplicative. Then for $p$ a prime and $m\geq 2$, we have $g(p^m)=f(p^m)-f(p^{m-1})=\left(\frac{p}{p-1}\right)^4\frac{1}{p^m}$. For $m=1$, we find $g(p)=f(p)-1=\left(\frac{p}{p-1}\right)^4(1+1/p)-1$. Thus
\begin{align*}
\sum_{d\leq \xi^2}\frac{f(d)}{d}\chi_0(d)&=\sum_{k\leq\xi^2}\frac{g(k)\chi_0(k)}{k}\sum_{d\leq \xi^2/k}\frac{\chi_0(d)}{d}\leq \sum_{d\leq \xi^2}\frac{\chi_0(d)}d\sum_{k=1}^\infty\frac{g(k)}{k}
\\&\lesssim 2\varkappa \log T\frac{\phi(q)}{q}\prod_p\left(1+\frac{5p^5-6p^4+5p^2-4p+1}{(p-1)^5p(p+1)}\right).
\end{align*}

Now we turn to the error term in \eqref{eq-Im_Kv/v}. We have
\begin{align*}
\frac{d^{5/2}\sum_{\rho|d}\frac1\rho}{\phi(d)^{7/2}}\ll \frac{1}{d}\prod_{p|d}(1-1/p)^{-9/2}\ll \frac{1}{d}\prod_{p|d}(1+p^{-3/4}),
\end{align*}
and hence the said error term is
\[\ll \frac1{(\log T)^{7/6}}\sum_{d\leq T}\frac{1}{d}\prod_{p|d}(1+p^{-1/2})\ll \frac1{(\log T)^{7/6}}\sum_{d\leq T}\frac1d\sum_{\rho|d}\rho^{-1/2}\leq \frac1{(\log T)^{1/6}}\zeta(3/2).\]

Therefore for any $0<v<1/\log \xi$:
\begin{align*}
\Im\left(\frac{K(v)}v\right)&\lesssim 2\varkappa\frak C_5^2\frac{q}{\phi(q)}\prod_p\left(1+\frac{5p^5-6p^4+5p^2-4p+1}{(p-1)^5p(p+1)}\right)
\\\numberthis\label{petitv}&\qquad\times\left[\left(\frac12+2\varkappa\right)\frak C_6(v\log T)^2+2\frak C_4\frak C_6(v\log T)\sqrt\varkappa\right].
\end{align*}
For commodity we define 
\[\frak C_7(v\log T):=\left(\frac12+2\varkappa\right)\frak C_6(v\log T)^2+2\frak C_4\frak C_6(v\log T)\sqrt\varkappa.\]

Now, we wish to study the quantity $\Im(K(v)/v)$ when $v$ is larger, namely $1/\log \xi <v< H$. We fix $d\leq \xi^2$ co-prime to $q$. Using \nameref{lemm12} with $r=2$ and using \eqref{lien-molli}, one finds that for $\delta|d^\infty$, $\delta\leq \xi$:
\begin{align*}
&\sum_{(n,qd)=1}\frac{\M(n\delta)\tau_{-1/2}(n\delta)}{n^{1+iv}}
\\&=\frac{-\tau_{-1/2}(\delta)}{(1-\theta)2\sqrt\pi \sqrt v\log \xi}\prod_{p|qd}(1-p^{-(1+iv)})^{-1/2} \left[\Delta_2\left(v\log\frac\xi\delta\right)-\Delta_2\left(v\log^+\frac{\xi^\theta}\delta\right)\right]
\\&\numberthis\label{eqintero}\qquad +\O\left(\tau_{1/2}(\delta)\prod_{p|d}(1+p^{-3/4})(\log T)^{-2/3}\right).
\end{align*}
For commodity we denote the quantity between brackets by $\Delta(\xi,\delta,\theta,v)$.

Putting \eqref{eqintero} in \eqref{eqa'}, this leads to 
\begin{align*}
A'_d(v)&=\sum_{\substack{d|\delta\delta'\\\delta\delta'|d^\infty\\\delta,\delta'\leq \xi}}\frac1{\delta^{1+iv}\delta'}\left(\frac{-\tau_{-1/2}(\delta)}{(1-\theta)2\sqrt\pi \sqrt v\log \xi}\prod_{p|qd}(1-p^{-(1+iv)})^{-1/2}\Delta(\xi,\delta,\theta,v)\right)B_d(\delta')
\\&\qquad+\O\left(\prod_{p|d}(1+p^{-3/4})\sum_{\substack{d|\delta\delta'\\\delta\delta'|d^\infty}}\frac{\tau_{1/2}(\delta)}{\delta\delta'(\log T)^{2/3}}|B_d(\delta')|\right).
\end{align*}
By \eqref{defC'}, the error term is 
\[\ll\frac{1}{(\log T)^{7/6}}\sqrt{\frac {d}{\phi(d)}}\frac{1}{\phi(d)}\prod_{p|d}(1+p^{-3/4})\ll\frac{1}{(\log T)^{7/6}}\frac{1}d\prod_{p|d}(1+p^{-3/4})^2.\]
Using the bound in \nameref{lemm_esti_delta}, we deduce that $\Delta(\xi,\delta,\theta,v)=\O(\frak A)$, and hence
\begin{align*}
\numberthis\label{eqADcarré} A'_d(v)^2&=\frac{1}{4(1-\theta)^2\pi v(\log \xi)^2}\prod_{p|qd}(1-p^{-1-iv})^{-1}\left(\sum_{\substack{d|\delta\delta'\\\delta\delta'|d^\infty\\\delta,\delta'\leq \xi}}\frac{\tau_{-1/2}(\delta)B_d(\delta')\Delta(\xi,\delta,\theta,v)}{\delta^{1+iv}\delta'}\right)^2
\\&+\O\left(\frac{\frak A}{(\log T)^{13/6}}\frac{1}{d^2}\prod_{p|d}(1+p^{-3/4})^3+\frac{1}{(\log T)^{7/3}}\frac{1}{d^2}\prod_{p|d}(1+p^{-3/4})^4\right).
\end{align*}
Noticing that $\phi_{iv}(d)=d^{1+iv}\prod_{p|d}(1-p^{-1-iv})$ and using this last equality in \eqref{eqExpK}, we find
\begin{align*}
&\Im\left(\frac{K(v)}{v}\right)=\O(\frak A /(v(\log T)^{7/6}))
\\&+\frac{1}{4(1-\theta)^2\pi (v\log \xi)^2}\sum_{\substack{d\leq \xi^2\\(d,q)=1}}\Im\left[\prod_{p|q}(1-p^{-1-iv})^{-1}\mathcal T^{iv/2}d^{1+iv}\left(\sum_{\substack{d|\delta\delta'\\\delta\delta'|d^\infty\\\delta,\delta'\leq \xi}}\frac{\tau_{-1/2}(\delta)B_d(\delta')\Delta(\xi,\delta,\theta,v)}{\delta^{1+iv}\delta'}\right)^2\right].
\end{align*}
Using this together with \eqref{petitv} and using the change of variable $u=v\log T$ in the following integrals, we find the following 
\begin{align*}
&\int_0^H(H-v)\Im\left(\frac{K(v)}v\right)\diff v
\\&\lesssim 2\varkappa\frak C_5^2\frac{q}{\phi(q)}\prod_p\left(1+\frac{5p^5-6p^4+5p^2-4p+1}{(p-1)^5p(p+1)}\right)\int_0^{1/\varkappa}\frac{(\frak A-u)}{(\log T)^2}\frak C_7(u)\diff u+\frac{1}{4(1-\theta)^2\pi (\log \xi)^2}
\\&\numberthis\label{intKV}\times\sum_{\substack{d\leq \xi^2\\(d,q)=1}}\int_{1/\log \xi}^H(H-v)\Im\left[\prod_{p|q}(1-p^{-1-iv})^{-1}\mathcal T^{iv/2}d^{1+iv}\left(\sum_{\substack{d|\delta\delta'\\\delta\delta'|d^\infty\\\delta,\delta'\leq \xi}}\frac{\tau_{-1/2}(\delta)B_d(\delta')\Delta(\xi,\delta,\theta,v)/v}{\delta^{1+iv}\delta'}\right)^2\right]\diff v
\end{align*}
Finally we have, for $d\leq \xi^2$, $(q,d)=1$:
\begin{align*}
&\left|\int_{1/\log \xi}^H(H-v)\Im\left[\prod_{p|q}(1-p^{-1-iv})^{-1}\mathcal T^{iv/2}d^{1+iv}\left(\sum_{\substack{d|\delta\delta'\\\delta\delta'|d^\infty\\\delta,\delta'\leq \xi}}\frac{\tau_{-1/2}(\delta)B_d(\delta')\Delta(\xi,\delta,\theta,v)/v}{\delta^{1+iv}\delta'}\right)^2\right]\diff v\right|
\\\leq & d\left|\int_{1/\log \xi}^H(H-v)(\sqrt{\mathcal T}d)^{iv}\prod_{p|q}(1-p^{-1-iv})^{-1}\left(\sum_{\substack{d|\delta\delta'\\\delta\delta'|d^\infty\\\delta,\delta'\leq \xi}}\frac{\tau_{-1/2}(\delta)B_d(\delta')\Delta(\xi,\delta,\theta,v)/v}{\delta^{1+iv}\delta'}\right)^2\diff v\right|
\\\leq & d\frac{q}{\phi(q)}\sum_{\substack{d|\delta_1\delta_1'\\\delta_1\delta_1'|d^\infty\\\delta_1,\delta_1'\leq \xi}}\sum_{\substack{d|\delta_2\delta_2'\\\delta_2\delta_2'|d^\infty\\\delta_2,\delta_2'\leq \xi}}\frac{\tau_{1/2}(\delta_1)\tau_{1/2}(\delta_2)}{\delta_1\delta_1'\delta_2\delta_2'}|B_d(\delta_1')B_d(\delta_2')|
\\\numberthis\label{lastone}&\qquad\times\left|\int_{1/\log \xi}^H(H-v)\left(\frac{\sqrt{\mathcal T} d}{\delta_1\delta_2}\right)^{iv}\frac{\Delta(\xi,\delta_1,\theta,v)}v\frac{\Delta(\xi,\delta_2,\theta,v)}v\diff v\right|.
\end{align*}
Using \nameref{lemm_esti_delta}, we find that\footnote{Note that this is precisely here that the quality of our lower bound for the proportion of critical zeros is at stake. Selberg's paper led to an error term $\ll\sqrt{v\log \xi}$, and hence the integral of our \eqref{gainSel} would become $\int_{1/\log \xi}^H (H-v)/\sqrt{v}\diff v$. This would produce a term $\asymp \frak A^{3/2}$, while our biggest term here is $\asymp\frak A\log\frak A$.}
\[\Delta(\xi,\delta,\theta,v)=\left\{\begin{array}{ll}
2\sqrt\pi e^{-3i\pi/4}v\log\frac{\xi}\delta+\pb\left(\frac{32}{3}\right) & \text{if }\xi\geq \delta\geq \xi^\theta,
\\ 2\sqrt\pi e^{-3i\pi/4}v(1-\theta)\log\xi+\pb\left(\frac{32}{3}\right)& \text{if }1\leq \delta<\xi^\theta.
\end{array}\right.\]

Now we split the sums of \eqref{lastone} in four, according to whether $\delta_i\leq \xi^\theta$, $i=1,2$, or not. In the case $\xi\geq \delta_1,\delta_2>\xi^\theta$, the \eqref{lastone} integral term is equal to
\begin{align*}
&\left|\int_{1/\log \xi}^H\int_{1/\log \xi}^u\left(\frac{\sqrt{\mathcal T} d}{\delta_1\delta_2}\right)^{iv}\left(4\pi i v^2 \log\frac{\xi}{\delta_1}\log\frac{\xi}{\delta_2} +\pb\left(\frac{128}3\sqrt\pi (1-\theta)v\log\xi+\frac{1024}9\right)\right)\frac{\diff v}{v^2}\diff u\right|
\\&\lesssim\frac{4\pi (1-\theta)^2\log^2\xi}{\log(\sqrt{\mathcal T} d/(\delta_1\delta_2))}\left|\int_{1/\log \xi}^H\left(\left(\frac{\sqrt{\mathcal T} d}{\delta_1\delta_2}\right)^{iu}-\left(\frac{\sqrt{\mathcal T} d}{\delta_1\delta_2}\right)^{i/\log \xi}\right)\diff u\right|+\frac{1024}9\int_{1/\log \xi}^H\frac{H-v}{v^2}\diff v
\\\numberthis\label{gainSel}&\qquad +\frac{128}3 \sqrt\pi (1-\theta)\log\xi \int_{1/\log\xi}^H\frac{H-v}v\diff v
\\&\leq 4\pi(1-\theta)^2\log^2\xi\left(\frac{H-1/\log \xi}{\log(T^{1/2-2\varkappa})}+\frac2{\log(T^{1/2-2\varkappa})^2}\right) +\frac{1024}9 (H\log\xi-\log(H\log \xi)-1)
\\&\qquad +\frac{128}3\sqrt\pi (1-\theta)\log\xi\left(\frac1{\log \xi}-H+H\log(H\log \xi)\right)
\\&= 4\pi(1-\theta)^2\varkappa^2\left(\frac{\frak A-1/\varkappa}{1/2-2\varkappa}+\frac{2}{(1/2-2\varkappa)^2}\right)+\frac{1024}9(\frak A\varkappa-\log(\frak A\varkappa)-1)
\\&\qquad +\frac{128}3\sqrt\pi(1-\theta)(1-\frak A\varkappa+\frak A\varkappa\log(\frak A\varkappa)).
\end{align*}
The three other cases lead to the same upper bound. Thus, summing \eqref{lastone} over $d\leq \xi^2$, $(q,d)=1$, we find that the second summand of \eqref{intKV} is 
\begin{align*}
&\lesssim\frac{1}{(\log T)^2}\Biggl\{\frac{\frak A-1/\varkappa}{1/2-2\varkappa}+\frac{2}{(1/2-2\varkappa)^2}+\frac{256/9}{(1-\theta)^2\pi\varkappa^2}(\frak A\varkappa-\log(\frak A\varkappa)-1)
\\&\qquad+\frac{32/3}{\sqrt\pi (1-\theta)\varkappa^2}(1-\frak A\varkappa+\frak A\varkappa\log(\frak A\varkappa))\Biggr\}\frak C_5^2\frac{q^2}{\phi(q)^2}\sum_{d\leq \xi^2}\frac{d^2\chi_0(d)}{\phi(d)\log T}\left(\sum_{\substack{d|\delta\delta'\\\delta\delta'|d^\infty}}\frac{\tau_{1/2}(\delta)\tau_{1/2}(\delta')}{\delta\delta'}\right)^2
\\&\lesssim\frac{1}{(\log T)^2}\Biggl\{\frac{\frak A-1/\varkappa}{1/2-2\varkappa}+\frac{2}{(1/2-2\varkappa)^2}+\frac{256/9}{(1-\theta)^2\pi\varkappa^2}(\frak A\varkappa-\log(\frak A\varkappa)-1)
\\&\qquad+\frac{32/3}{\sqrt\pi (1-\theta)\varkappa^2}(1-\frak A\varkappa+\frak A\varkappa\log(\frak A\varkappa))\Biggr\}\frak C_5^2\frac{q}{\phi(q)} \times 2\varkappa \prod_{p}\left(1+\frac{3p^2-3p+1}{p^4-3p^3+3p^2-p}\right).
\end{align*}
Here we used \eqref{defC'} to bound $B_d(\delta)$ and \eqref{eqprodeul} to estimate the sum over $d$, which is shown to be $\leq\sum_{d\leq \xi^2}d^2\chi_0(d)/\phi(d)^3$. Using \eqref{eqIntfacile} and \eqref{intKV}, we find the expected result : 
\[\int_T^{2T}|I(t,H)|^2\diff t\lesssim \frac{T}{(\log T)^2}8\frak C_5^2\left(\frak K_1 \frak A\log\frak A+\frak K_2\frak A+\frak K_3\log \frak A+\frak K_4\right),\]
where
\[\numberthis\label{deffrakK}\frak K_1:=\prod_{p}\left(1+\frac{3p^2-3p+1}{p^4-3p^3+3p^2-p}\right)\frac{32/3}{\sqrt\pi (1-\theta)},\]
\begin{align*}
\frak K_2&:=\varkappa\prod_p\left(1+\frac{5p^5-6p^4+5p^2-4p+1}{(p-1)^5p(p+1)}\right)\int_0^{1/\varkappa}\frak C_7(v)\diff v
\\&+  \prod_{p}\left(1+\frac{3p^2-3p+1}{p^4-3p^3+3p^2-p}\right)\left[\frac{\varkappa}{1/2-2\varkappa}+\frac{256/9}{(1-\theta)^2\pi}+\frac{32/3}{\sqrt\pi (1-\theta)}(\log\varkappa-\varkappa)\right],
\end{align*}
\[\frak K_3:=- \prod_{p}\left(1+\frac{3p^2-3p+1}{p^4-3p^3+3p^2-p}\right)\frac{256/9}{(1-\theta)^2\pi\varkappa},\]
and finally
\begin{align*}
\frak K_4&:= -\varkappa\prod_p\left(1+\frac{5p^5-6p^4+5p^2-4p+1}{(p-1)^5p(p+1)}\right)\int_0^{1/\varkappa}v\frak C_7(v)\diff v
\\&+\prod_{p}\left(1+\frac{3p^2-3p+1}{p^4-3p^3+3p^2-p}\right)\left[\frac{4\varkappa-1/2}{(1/2-2\varkappa)^2}-\frac{256/9}{(1-\theta)^2\pi\varkappa}(1+\log\varkappa)+\frac{32/3}{\sqrt\pi (1-\theta)\varkappa}\right].
\end{align*}
This concludes.
\end{proo}
\section{Computations}
Now it only remains to explain how we did the computations of \nameref{application}.
\begin{proo}{of \nameref{application}}
In the case of small fixed $N$, we have to compute $\frak C_1$. This is easily done by noticing that $\frak C_7(v)$ is an increasing and positive function of $v$, for any fixed $\theta\in (0,1)$, as standard arguments show. Therefore, it is easy to bound the integrals of $\frak K_2$ and $\frak K_4$ by the rectangle method, for example. We chose to bound these integrals by using $100$ regular rectangles. Then, the choice of $\frak A$ was done by Sage\footnote{Code available \href{https://github.com/JeremyDous/Constante-Selberg/blob/f2fd5248fd063eae0c6857e22c2ff3a1e545ff0a/Code\%20Constante\%20Selberg.py}{here}.}, by finding a positive root of the derivative of\footnote{Except in the case $N=1$, where we used \eqref{borne_meilleure} instead.}
\[2\pi\left(\frac{1}{2\frak A}-4N\frac{\frak C_1(\frak A)+\frak C_2 }{\frak A^3}\right).\]
This gives a result for a fixed parameter $\theta$. To find a good $\theta$ we split the interval $(0,1)$ in $10^4$, and then applied the process to each of these $\theta$. The program then returned the one giving the best result.

In the case of large $N$, we have for some constant $C$
\[\kappa_F\approx 2\pi\left(\frac{1}{2\frak A}-4N C\frac{\frak A\log\frak A}{\frak A^3}\right),\]
which is optimized when $1/\frak A\approx N\frak A\log\frak A/\frak A^3$, and hence when $\frak A$ is of the form $\lambda N\log N$ for some positive number $\lambda$. Thus, we may say that 
\[\frak C_1(\frak A)\approx 8\frak C_5^2\frak K_1 \frak A\log \frak A\approx 8\frak C_5^2 \frak K_1 \lambda N\log(N)^2.\]
Because of the definition of these constants, we should take $\theta$ as small as possible, but this will makes $\frak C_2$ become very large. With this in mind, if $N$ is so large that $\frak C_2$ becomes negligible, we may write that heuristically
\[2\pi\left(\frac{1}{2\frak A}-4N\frac{\frak C_1(\frak A)+\frak C_2 }{\frak A^3}\right)\approx \frac{2\pi}{N\log N}\left(\frac{1}{2\lambda}-\frac{32\frak C_5^2\frak K_1}{\lambda^2} \right),\]
which is optimized by taking $\lambda=128\frak C_5^2\frak K_1=128\frak C_5^2(\theta)\frak K_1(\theta)$. Now, we simply inject this in the lower bound for $\kappa_F$ and find how large $N$ should be to compensate the size of $\frak C_2(\theta)$. 

We put $\frak A=\lambda(\theta) N\log N$, fix $\varepsilon>0$ and we take $0<\theta<1$ so that $1/(1-\theta)=1+\varepsilon$.

Observe that $(e^{\varrho_\theta}+e^{\varrho_\theta\theta})/\sqrt{\varrho_\theta}$ is an increasing function of $\theta$, and note that for $\varepsilon<1/3$, we have $\theta<1/4$. Thus
\[(1+\varepsilon)\frak C_5^-:=(1+\varepsilon)\frac{e^{\varrho_0}+1}{2\sqrt{\pi\varkappa\varrho_0}}\frac{\Gamma(1/4)}{\Gamma(3/4)}\leq\frak C_5(\theta)\leq (1+\varepsilon)\frac{e^{\varrho_{1/4}}+e^{\varrho_{1/4}/4}}{2\sqrt{\pi\varkappa\varrho_{1/4}}}\frac{\Gamma(1/4)}{\Gamma(3/4)}=:(1+\varepsilon)\frak C_5^+.\]
Moreover, $\frak K_1(\theta)=(1+\varepsilon)\frak K_1(0)$. Thus we have that 
\[(1+\varepsilon)^3\lambda^-:=128(\frak C_5^-)^2\frak K_1(0)(1+\varepsilon)^3\leq\lambda(\theta)\leq 128(\frak C_5^+)^2\frak K_1(0)(1+\varepsilon)^3=:(1+\varepsilon)^3\lambda^+.\]
We also have
\begin{align*}
\frak C_3(\theta)&\leq \left(\frac{1}{8\varkappa}+\frac32\right)\left(\frac{(e^{\varrho(1/4)}+e^{\varrho(1/4)/4})\Gamma(1/4)}{\sqrt{\pi\varrho(1/4)}\Gamma(3/4)}\right)^4(1+\varepsilon)^4\prod_p\left(1+\frac{3p^2-3p+1}{p^4-3p^3+3p^2-p}\right)
\\&=:(1+\varepsilon)^4\frak C_3^+.
\end{align*}
This implies that 
\[\frak C_2(\theta)\leq \frac{6(\frak C_3^++1+2\sqrt{\frak C_3^+})}{\varkappa^2}\frac{(1+\varepsilon)^6}{\varepsilon^2}=:\frac{(1+\varepsilon)^6}{\varepsilon^2}\frak C_2^+.\]

We also take 
\[\numberthis\label{lowboundN}N\geq \frac{\frak C_2^+}{(\lambda^-)^3\varepsilon^3}=:\frac{N_0}{\varepsilon^3}\]
so that $\frac{\frak C_2(\theta)}{\lambda(\theta)^3 N\log^2N}\leq \frac{\varepsilon}{\log^2N}$. Remark that $\frak K_3(\theta)\leq 0$ and
\[\frak K_4(\theta)\leq \prod_{p}\left(1+\frac{3p^2-3p+1}{p^4-3p^3+3p^2-p}\right)\left[\frac{256/9}{\pi\varkappa}(\log(1/\varkappa)-1)+\frac{32/3}{\sqrt\pi\varkappa}\right](1+\varepsilon)^2=:(1+\varepsilon)^2\frak K_4^+.\]
Finally, note that for $0\leq v\leq 8$, $\theta\in(0,1)$, we have $1/2\leq \rho(v,\theta)\leq 1$ and hence
\[\frak C_6(v,\theta)\leq (1+\varepsilon)\frac{e}{\sqrt{\pi\varkappa/2}}\left(\sqrt{2v\pi\varkappa}+\frac{\Gamma(1/4)}{\Gamma(3/4)}\right)=:(1+\varepsilon)\frak C_6^+(v).\]
Since $\frak C_4(\theta)\leq 13/(5\sqrt\pi)(1+\varepsilon)^2=:(1+\varepsilon)^2\frak C_4^+$, we deduce that 
\begin{align*}
&\frak K_2(\theta)\leq \Biggl\{\varkappa\prod_p\left(1+\frac{5p^5-6p^4+5p^2-4p+1}{(p-1)^5p(p+1)}\right)\int_0^{1/\varkappa}\left(\frak C_6^+(v)^2\left(\frac{1}2+2\varkappa\right)+2\frak C_4^+\frak C_6^+(v)\sqrt\varkappa\right) \diff v
\\&+\prod_{p}\left(1+\frac{3p^2-3p+1}{p^4-3p^3+3p^2-p}\right)\left[\frac{\varkappa}{1/2-2\varkappa}+\frac{256}{9\pi}\right]\Biggr\}(1+\varepsilon)^3=:(1+\varepsilon)^3\frak K_2^+,
\end{align*}
Thus, for $N$ satisfying \eqref{lowboundN}, we have
\begin{align*}
4\frac{\frak C_1(\frak A,\theta)+\frak C_2(\theta)}{\lambda(\theta)^3N\log^2N}&\leq \frac{1}{4\lambda(\theta)}+\frac{\log\log N}{4\lambda^-\log N}+\frac{\log \lambda^++3\varepsilon+\frak K_2^+/\frak K_1(0)}{4\lambda^-\log N}+\frac{\frak K_4^+/\frak K_1(0)}{4(\lambda^-)^2N\log^2N}+\frac{4\varepsilon}{\log ^2N}.
\end{align*}
Therefore, for any such $N$, we find that 
\begin{align*}
\kappa_F&\geq \frac{2\pi}{N\log N}\left(\frac{1}{2\lambda(\theta)}-4\frac{\frak C_1(\frak A,\theta)+\frak C_2(\theta)}{\lambda(\theta)^3N\log^2N}\right)
\\&\geq \frac{2\pi}{N\log N}\left(\frac{1}{4\lambda^+(1+\varepsilon)^3}-\frac{\log\log N}{4\lambda^-\log N}-\frac{\log \lambda^++1+\frak K_2^+/\frak K_1(0)}{4\lambda^-\log N}-\frac{\frak K_4^+/\frak K_1(0)}{4(\lambda^-)^2N\log^2N}-\frac{4\varepsilon}{\log ^2N}\right).
\end{align*}
Using the fact that $1/(1+\varepsilon)^3\geq 1-3\varepsilon$, a computation then leads to the expected result.
\end{proo}
\bibliographystyle{abbrv}
\bibliography{biblio}
\labo
\end{document}